\documentclass[12pt]{amsart}
\usepackage{amssymb,amsfonts,graphicx,mathrsfs, amsmath}
\usepackage{color}

\definecolor{myblue}{rgb}{0,0,0.6}
\definecolor{myred}{rgb}{0.9,0,0.1}

\def\foorp{{\hfill$\spadesuit$}}

\def\inv{{^{-1}}}
\def\RR{{\mathbb R}}
\def\PP{{\mathbb P}}

\def\CC{\mathbb C}
\def\HH{\mathbb H}
\def\ZZ{\mathbb Z}
\def\MM{\mathbb M}

\def\HS{\mathcal{H}\mathcal{S}}

\def\be{\begin{equation}}
\def\beq#1{\begin{equation}\label{#1}}
\def\ee{\end{equation}}
\def\bea{\begin{eqnarray}}
\def\beqa#1{\begin{eqnarray}\label{#1}}
\def\eea{\end{eqnarray}}
\def\ba{\begin{array}}
\def\ea{\end{array}}

\def\cA{{\mathcal A}}
\def\cB{{\mathcal B}}
\def\cC{{\mathcal C}}

\def\cE{{\mathcal E}}
\def\cF{{\mathcal F}}

\def\cH{{\mathcal H}}
\def\cJ{{\mathcal J}}

\def\cM{{\mathcal M}}

\def\cS{{\mathcal S}}
\def\cU{{\mathcal U}}
\def\cV{{\mathcal V}}

\def\bds{_{-\infty}^\infty}
\def\supp{\hbox{supp}}

\def\bm{{\mathbf m}}

\renewcommand\cV{\mathscr{V}}
\renewcommand\cS{\mathscr{S}}
\newcommand\cHH{\mathscr{H}}
\newcommand\cBB{\mathscr{B}}

\renewcommand\cM{\mathscr{M}}

\newtheorem{theorem}{Theorem}
\newtheorem{proposition}{Proposition}
\newtheorem{lemma}{Lemma}
\newtheorem{corollary}{Corollary}
\newtheorem{Def}{Definition}
\theoremstyle{remark}
\newtheorem{remark}{Remark}
\newtheorem{example}{Example}
\begin{document}

\title[Operators in time-frequency domain]{Representation  of operators in  the time-frequency domain and generalized {G}abor multipliers }


\author{Monika D\"orfler}
\address{Acoustics Research Institute, Austrian Academy of Science, Wohllebengasse 12-14, A-1040 Vienna,
Austria}
\email{monid@kfs.oeaw.ac.at}
\thanks{The first author has been supported by the WWTF project MA07-025.}

\author{Bruno Torr\'esani}
\address{Laboratoire d'Analyse, Topologie
et Probabilit\'es, Centre de Math\'ematique et d'Infor\-matique,
39 rue Joliot-Curie, 13453 Marseille cedex 13, France}
\curraddr{}
\email{bruno.torresani@univ-provence.fr}
\thanks{}



\date{\today}

\dedicatory{}

\begin{abstract}
 Starting from a general operator representation in the time-frequency domain, this paper addresses the problem of approximating linear operators by operators that are
diagonal or band-diagonal with respect to  Gabor frames. A characterization of operators that can be realized as
Gabor multipliers is given  and necessary conditions for
the existence of (Hilbert-Schmidt) optimal Gabor multiplier
approximations are discussed and an efficient method for the calculation of an operator's best approximation by a Gabor multiplier is derived. The spreading function of Gabor multipliers yields new error estimates for these approximations.  Generalizations (multiple Gabor multipliers) 
are introduced for better approximation of overspread operators.   The Riesz property of the projection operators involved in generalized Gabor multipliers is characterized, and a method for obtaining an operator's best approximation by a multiple Gabor multiplier is suggested. Finally, it is shown that in certain situations, generalized Gabor multipliers reduce to a finite sum of regular Gabor multipliers with adapted windows.    
\end{abstract}

\maketitle

\section{Introduction}\label{se:intro}
The goal of time-frequency analysis is to provide efficient
representations for functions or distributions in terms of
decompositions such as
\[{
 f =
\sum_{\lambda\in\Lambda} \langle f,g_\lambda \rangle h_\lambda\ .
}\]
Here, $f$ is expanded as a weighted sum of {\em atoms} $ h_\lambda$
well localized in
both time and frequency domains. The time-frequency coefficients
$ \langle f,g_ \lambda \rangle$ characterize the function
under investigation, and a {\em synthesis map} usually allows the
reconstruction of the original function $f$.\\
Concrete applications can be found mostly in signal analysis and processing
(see~\cite{Carmona98practical,Daubechies92ten,Mallat98wavelet} and
references therein), but recent works
in different areas such as numerical analysis may also be mentioned
(see for example~\cite{fest98,fest03} and
references therein).

Time-frequency analysis of operators, originating in the work on communication channels of Bello~\cite{be63}, Kailath~\cite{ka62} and Zadeh~\cite{za61}, has
enjoyed increasing interest during the last few years,~\cite{hlma02,kopf06,pfwa05-1,baporive05}. Efficient
time-frequency operator representation is a  challenging task, and often
the intuitively appealing approach of operator approximation 
by modification of the time-frequency coefficients before reconstruction
is the method of choice.  If the modification of the coefficients is
confined to be multiplicative, this approach leads to the model of
time-frequency multipliers, as discussed in Section~\ref{se:mult}.
The class of operators that may be well represented by
time-frequency multipliers depends on the choice of the
parameters involved and is  restricted to operators performing
only small time-shifts or modulations. 

The work in this paper is  inspired by a general operator
representation in the time-frequency domain via a
\emph{twisted convolution}. It turns out, that this
representation, respecting the underlying structure of the
Heisenberg group, has an interesting connection to the
so-called spreading function representation of operators. An
operator's spreading function comprises the amount of time-shifts and
modulations, i.e. of time-frequency-shifts, effected by the operator.
Its investigation is hence decisive in the study of time-frequency
multipliers and their generalizations. Although no direct
discretization of the continuous representation by an operator's
spreading function is possible, the twisted convolution turns out
to play an important role in the generalizations of time-frequency
multipliers.  In the main section of this article, we introduce a
general model for {\emph multiple Gabor multipliers (MGM)}, which uses several  synthesis windows simultaneously. Thus, by jointly
adapting the respective masks, more general operators may be
well-represented than by regular Gabor multipliers. Specifying to
a separable mask in the modification of time-frequency coefficients
within  MGM, as well as a specific sampling lattice for the
synthesis windows, it turns out, that the MGM reduces to one or the sum of a 
finite number of regular Gabor multipliers with  adapted synthesis windows. \\
For the sake of generality, most statements are given in a
\emph{ Gelfand-triple}, rather than a pure Hilbert space setting.
This choice bears several advantages. First of all, many important
operators and signals may not be described in a Hilbert-space
setting, starting from simple operators as the identity.
Furthermore, by using distributions, continuous and discrete
concepts may be considered together. Finally, the Gelfand-triple
setting often allows for short-cut proofs of statements formulated
in a general context. 

This paper is organized as follows.  The next section gives a review of
the time-frequency plane and the corresponding continuous and discrete
transforms. We then introduce the concept of Gelfand triples, which
will allow us to consider operators beyond the Hilbert-Schmidt
framework. The section closes with the important statement on
operator-representation in the time-frequency domain via twisted
convolution with an operator's spreading function. Section~\ref{Se:tfmult}
introduces time-frequency multipliers and gives a criterion for
their ability to approximate linear operators. A fast method for the calculation of an operator's best approximation
by a Gabor multiplier in Hilbert-Schmidt sense is suggested. 
 Section~3 introduces
generalizations of Gabor multipliers. The operators in the construction
of MGM are investigated and a criterion for their Riesz basis property
in the space of Hilbert-Schmidt operators is given. We mention some
connections to classical Gabor frames. A numerical example
concludes the discussion of general MGM. In the final section,
TST (twisted spline type) spreading functions are introduced. It
is shown, that under certain conditions, a MGM reduces to a regular
Gabor multiplier with an adapted window or a finite sum of regular
multipliers with the same mask and adapted windows.

\section{Operators from the Time-frequency point of view}\label{se:TFop}
Whenever one is interested in time-localized frequency information in
a signal or operator, one is naturally led to the notion of the
time-frequeny plane, which, in turn, is closely related to the
Weyl-Heisenberg group.
\subsection{Preliminaries: the time-frequency plane}
The starting point of our operator analysis is the so-called
{\em spreading function operator representation}. This operator
representation expresses linear operators as  a sum (in 
a sense to be specified below) of time-frequency shifts
$\pi(b,\nu) = M_\nu T_b$. Here,  the translation and modulation
operators are defined as
$$
T_b f(t) = f(t-b)\ , \quad M_\nu f(t) = e^{2i\pi\nu t}f(t)\ ,
\quad f\in \mathbf{L}^2(\RR)\ .
$$
These (unitary) operators generate a group, called the Weyl-Heisenberg group
\begin{equation}
\label{fo:WH.group}
\HH = \left\{(b,\nu,\varphi)\in\RR\times\RR\times [0,1[\right\}\ ,
\end{equation}
with group multiplication
\begin{equation}
(b,\nu,\varphi)(b',\nu',\varphi') = (b+b',\nu+\nu',\varphi+\varphi'-\nu'b)\ .
\end{equation}
The specific quotient space $\PP= \HH / [0,1]$ of the Weyl-Heisenberg
group is called {\em phase space}, or {\em time-frequency plane},
which plays  a central role in the subsequent analysis.  Details 
on the Weyl-Heisenberg group and the time-frequency plane may be found
in~\cite{Folland89harmonic,Schempp86harmonic}. In the current
article, we shall limit ourselves to the basic irreducible unitary
representation of $\HH$ on $\mathbf{L}^2(\RR)$, denoted by $\pi^o$, and
defined by
\begin{equation}
\pi^o(b,\nu,\varphi) = e^{2i\pi\varphi}M_\nu T_b\ .
\end{equation}
 By $\pi(b,\nu) = \pi^o(b,\nu,0)$  we  denote the restriction to
the phase space.
We refer to~\cite{Dorfler07spreading}
or~\cite[Chapter~9]{gr01} for a more detailed
analysis of this quotient operation.

\medskip
The left-regular (and right-regular) representation generally
plays a central role in group representation theory. By unimodularity of the 
Weyl-Heisenberg group, its left and right
regular representations coincide. We thus focus on the left-regular one,
acting on $\mathbf{L}^2(\HH)$ and defined by
\begin{equation}
\label{fo:left.regular}
\big[L(b',\nu',\varphi')F\big](b,\nu, \varphi) = 
F(b-b',\nu-\nu',\varphi-\varphi' +b'(\nu-\nu'))\ .
\end{equation}
Denote by $\mu$ the Haar measure.
Given $F,G\in \mathbf{L}^2(\HH,d\mu)$, the associated (left) convolution
product is the bounded function $F*G$, given by
\[
(F * G)(b,\nu,\varphi) = \int_\HH F(h)
\big[L(b,\nu,\varphi)G\big](h)\, d\mu(h).\] After quotienting out the phase term, this yields the
{\em twisted convolution} on $\mathbf{L}^2(\PP)$:
\begin{equation}
\label{fo:twisted.conv}
(F\natural G) (b,\nu) = \int\bds\int\bds
F(b',\nu') G(b-b',\nu-\nu') e^{-2i\pi b'(\nu-\nu')}\,db'd\nu'\ .
\end{equation}
The twisted convolution, which admits a nice interpretation in terms of
group Plancherel theory~\cite{Dorfler07spreading} is
non-commutative (which reflects the non-Abelianess of $\HH$)
but associative. It satisfies the usual Young inequalities,
but is in some sense nicer than the usual convolution,
since $\mathbf{L}^2(\RR^2) \natural \mathbf{L}^2(\RR^2) \subset
\mathbf{L}^2(\RR^2)$ (see~\cite{Folland89harmonic} for details).

\medskip
As explained in~\cite{Grossmann85transforms,Grossmann86transforms}
(see also~\cite{Fuhr05abstract} for a review), the representation
$\pi^o$ is unitarily equivalent to a subrepresentation of the left
regular representation. The representation coefficient is given  by a
variant of the {\em short time Fourier transform} (STFT), which we
define next.
\begin{Def}
Let $g\in\mathbf{L}^2(\RR)$, $g\ne 0$. The STFT of any $f\in\mathbf{L}^2(\RR)$
is the function on the phase space $\PP$ defined by
\begin{equation}
\cV_g f(b,\nu) = \langle f,\pi(b,\nu)g\rangle = \int\bds f(t) \overline{g}(t-b)
e^{-2i\pi\nu t}\,dt\ .
\end{equation}\end{Def}
This STFT is  obtained by quotienting out $[0,1]$ in the group
transform
\begin{equation}
\cV_g^o f(b,\nu,\varphi) = \langle f,\pi^o(b,\nu,\varphi)g\rangle\ .
\end{equation}
The integral transform $\cV_g^o$ intertwines $L$ and $\pi^o$,
i.e. $L(h) \cV_g^o = \cV_g^o\pi^o(h)$ for all $h\in\HH$. The latter
relation still holds true (up to a phase factor) when $\pi^o$ and
$\cV_g^o$ are replaced with $\pi$ and $\cV_g$ respectively.

It follows from the general theory of square-integrable
representations that for any $g\in\mathbf{L}^2(\RR)$, $g\ne 0$, the
transform $\cV_g^o$ is (a multiple of) an isometry $\mathbf{L}^2(\RR)\to
\mathbf{L}^2(\PP)$, and thus left invertible by the adjoint transform (up to a
constant factor). More precisely, given $h\in\mathbf{L}^2(\RR)$ such that
$\langle g,h\rangle\ne 0$, one has for all $f\in\mathbf{L}^2(\RR)$
\begin{equation}
\label{fo:STFT.inv}
f = \frac1{\langle h,g\rangle}\,
\int_\PP \cV_gf(b,\nu)\,\pi(b,\nu) h\,dbd\nu\ .
\end{equation}
We refer to~\cite{Carmona98practical,gr01} for more
details on the STFT and signal processing applications.

\medskip
The STFT, being a continuous transform,  is not well adapted for
numerical calculations, and, for practical issues, is  replaced by the Gabor
transform, which is a sampled version of it. To fix notation, we
outline  some steps of the Gabor frame theory and refer
to~\cite{Daubechies92ten,gr01} for a detailed
account.
\begin{Def}[Gabor transform]\label{Def:GabT}
Given $g\in\mathbf{L}^2(\RR)$  and two constants $b_0,\nu_0\in\RR^+$,
the corresponding Gabor transform associates with any
$f\in\mathbf{L}^2(\RR)$ the sequence of Gabor coefficients
\begin{equation}
\cV_gf(mb_0,n\nu_0) = \langle f,M_{n\nu_0} T_{mb_0} g\rangle
=  \langle f,g_{mn}\rangle\ ,
\end{equation}
where the functions
$
g_{mn}=M_{n\nu_0} T_{mb_0} g
$
are  the Gabor atoms associated to $g$ and the lattice
constants $b_0,\nu_0$.\end{Def}

Whenever the Gabor  atoms associated to $g$ and the given lattice $\Lambda = b_0 \mathbb{Z}\times \nu_0\mathbb{Z}$ form
a frame,\footnote{The operator \[S_{g }f  =
  \sum_{m,n\in\mathbb{Z}}\langle f, M_{mb_0}T_{n\nu_0}g\rangle
  M_{mb_0}T_{n\nu_0}g \] is the {\it frame operator} corresponding to
  $g$ and the lattice defined by $(b_0, \nu_0)$. If $S_g $ is
  invertible on $L^2 (\mathbb{R} )$, the family of time-frequency shifted atoms
  $M_{mb_0}T_{n\nu_0}g$, $m,n\in\mathbb{Z}$,  is a \emph{Gabor frame}
  for  $L^2 (\mathbb{R} )$.}
the Gabor transform is left invertible, and there exists
$h\in\mathbf{L}^2(\RR)$ such that any $f\in\mathbf{L}^2(\RR)$
may be expanded as
\begin{equation}
f = \sum_{m,n} \cV_gf(mb_0,n\nu_0)h_{mn}\ .
\end{equation}

\subsection{The Gelfand triple $( S_0, \mathbf{L}^2, S_0')$}
\label{Sec: Gelftrip}
We next set up a framework for the exact description of operators we
are interested in. In fact, by their property of being compact
operators, the Hilbert space of Hilbert-Schmidt operators turns out to
be far too restrictive to contain most operators of practical
interest, starting from the identity.
Although the classical triple $(\cS, \mathbf{L}^2, \cS')$ might seem
to be the appropriate choice  of generalization, we prefer to resort
to the Gelfand triple $( S_0, \mathbf{L}^2, S_0')$, which has proved
to be more adapted to a time-frequency environment. Additionally, the
Banach space property of $S_0$ guarantees a technically less elaborate
account.
\begin{Def}[$S_0$]\label{defS0}
Let $\cS (\mathbb{R})$ denote the Schwartz class.
Fix a  non-zero ``window'' function $\varphi \in \cS (\mathbb{R}   )$.
The space $S_0 (\mathbb{R} ) $ is given by
\[S_0 (\mathbb{R} ) = \{ f\in \mathbf{L}^2(\mathbb{R})  : 
\| f\|_{S_0}:= \|\mathcal{V}_{\varphi}
f\|_{L^1(\mathbb{R}^2)}<\infty\}.\]
\end{Def}

The following proposition summarizes some properties of
$S_0(\mathbb{R} )$ and its dual, the distribution space
$S_0'(\mathbb{R} )$.
\begin{proposition}\label{Pro:S0}
$S_0(\mathbb{R}   )$ is a Banach space and densely embedded in
  $\mathbf{L}^2(\mathbb{R})  $. The definition of
  $S_0(\mathbb{R})$ is  independent of the window $\varphi \in
  \cS(\mathbb{R} ) $, and different choices of $\varphi \in \cS(\RR) $
  yield equivalent norms on $S_0(\mathbb{R} )$. 

By duality, $\mathbf{L}^2(\mathbb{R})  $ is densely and
weak$^{\ast}$-continuously embedded in $S_0'(\mathbb{R}   )$ and can
also be characterized by  the norm $\| f\|_{S_0'} =
\|\mathcal{V}_{\varphi} f\|_{\mathbf{L}^{\infty}}$. 
\end{proposition}
In other words, the three spaces $(S_0(\mathbb{R} ),
\mathbf{L}^2(\mathbb{R} ) , S_0'(\mathbb{R} ))$ Represent   a 
special case of a \emph{Gelfand triple}~\cite{gelfand-vilenkin}
or Rigged Hilbert space.
For a proof, equivalent characterizations, and more results on $S_0$
we refer to ~\cite{feichtinger81,FZ98a,feichtinger-kozek98}.

Via an isomorphism between integral kernels in the Banach
spaces $S_0, S_0'$ and the operator spaces of bounded operators
$S_0'\mapsto S_0$ and $S_0\mapsto S_0'$, we obtain, together
with the Hilbert space  of Hilbert-Schmidt operators, a
Gelfand triple of operator spaces, as follows.
We denote by $\cBB$ the family of operators that are bounded
$S_0'\to S_0$ and by $\cBB'$ the family of operators that are bounded
$S_0\to S_0'$.
{
We have the following correspondence between these operator
classes and their integral kernels $\kappa$: }
\[
H\in (\cBB, \cHH, \cBB') \longleftrightarrow \kappa_H \in
(S_0(\mathbb{R} ), \mathbf{L}^2(\mathbb{R} ) , S_0'(\mathbb{R} ))\ .
\]
We will make  use of the principle of \emph{unitary Gelfand triple isomorphisms}, described for the Gelfand triples just introduced in~\cite{feichtinger-kozek98}. The basic idea is the extension of a  unitary isomorphism between $\mathbf{L}^2$-spaces to isomorphisms between the spaces forming the Gelfand triple. In fact, it may be shown, that it suffices to verify unitarity of a given isomorphic operator on the (dense) subspace $S_0$ in order to obtain a unitary Gelfand triple isomorphism, see~\cite[Corollary~7.3.4]{feichtinger-kozek98}. The most prominent examples for a unitary Gelfand triple isomorphism are the Fourier transform and the partial Fourier transform. 
For all  further details on the Gelfand triples just introduced,
we again refer to~\cite{feichtinger-kozek98},   
{
only mentioning here, that one important reason for investigating
operator representations on the level of Gelfand triples instead
of just a Hilbert  space framework is the fact, that $S_0'$
contains distributions such as the Dirac functionals, Shah
distributions or just pure frequencies and $\cBB'$ contains operators of great importance in signal processing, e.g.  convolution, the identity or just time-frequency shifts.

Subsequently, we will usually assume that the analysis and 
synthesis windows $g,h$ are in $S_0$. This is a rather mild condition,
which has almost become the canonical choice in Gabor analysis, for
many good reasons. Among others, this choice guarantees a beautiful
correspondence between the $\ell^p$-spaces and corresponding modulation space~\cite{gr01}. In the $\ell^2$-case this means, that 
the sequence of Gabor atoms generated from time frequency translates of an $S_0$ window
on an arbitrary lattice $\Lambda$  is automatically
a Bessel sequence (in such a case, the window is
termed ``Bessel atom''), which is not true for general
$\mathbf{L}^2$-windows.\\

The Banach spaces $S_0$ and $S_0'$ may also be interpreted as  Wiener amalgam spaces~\cite[Section~3.2.2]{FZ98a}. These time-frequency homogeneous spaces are defined as follows. Let $\cF\mathbf{L}^1$ denote the Fourier image of integrable functions and let a compact function $\phi\in\cF\mathbf{L}^1 (\RR )$ with $\sum_{n\in\ZZ} \phi(x-n ) \equiv 1$ be given. Then, for $\mathbf{ X} (\RR ) = \cF\mathbf{L}^1(\RR )$ or $\mathbf{ X }(\RR ) = \mathbf{C}(\RR )$, i.e. the space of continuous functions on $\RR $, or any of the Lesbesgue spaces,  we define, for $p\in [1,\infty )$, with the usual modification for $p=\infty$:
\begin{equation}\label{fo:amalgDef}
\mathbf{W}(\mathbf{ X} ,\ell^p ) = \big\{f \in	\mathbf{ X}_{loc}:\|f\|_{\mathbf{W}(\mathbf{ X} ,\ell^p )} = \big(\sum_{n\in\ZZ}\|f T_n\phi\|_{\mathbf{X}}^p\big)^{1/p}<\infty\big\}
\end{equation}
Now, $S_0 = \mathbf{W} (\cF\mathbf{L}^1,\ell^1 )$ and $S_0' = \mathbf{W} (\cF\mathbf{L}^{\infty},\ell^{\infty} )$, see~\cite[Section~3.2.2]{FZ98a}.\\ 

%
\subsection{The spreading function representation and its connections
to the STFT}\label{Se:sprep}
The so-called {\em spreading function} representation, closely related to the integrated Schr\"odinger representation~\cite[Section~9.2]{gr01},
expresses operators in $(\cBB, \cHH, \cBB')$ as a sum of
time-frequency shifts. More precisely, one
has (see~\cite[Chapter~9]{gr01}):
\begin{theorem}
\label{th:SFrep}
Let $H\in (\cBB,\cHH,\cBB')$; then there exists a spreading function
$\eta_H$ in  $(S_0(\RR^2),\mathbf{L}^2(\RR^2), S_0'(\RR^2))$ such that
\begin{equation}
\label{fo:SFrep}
H = \int\bds\int\bds \eta_H(b,\nu) \pi(b,\nu)\, dbd\nu\ .
\end{equation}
For $H\in\cH$, the correspondence $H\leftrightarrow\eta_H$ is isometric,
i.e. $\|H\|_\cHH = \|\eta_H\|_{\mathbf{L}^2(\PP)}$.
\end{theorem}
\begin{remark}\rm
{
For $H\in\cBB$, the decomposition given in \eqref{fo:SFrep} is
absolutely convergent, whereas, for $H\in\cBB'$, it holds in the
weak sense of bilinear forms on $S_0$.\\
}
When $\eta_H\in {\mathbf L}^2(\PP)$, $H$ is a Hilbert-Schmidt operator, and
the above integral is defined as a Bochner integral.
\end{remark}

The spreading function is intimately related to the integral kernel
$\kappa = \kappa_H$ of $H$ via
\begin{equation*}
\label{fo:SF2Ker2SF}
\hspace{-0.8cm}\eta_H(b,\nu)=\int\bds \kappa_H(t,t-b) e^{-2i\pi\nu t}\, dt
\mbox{ and }
\kappa_H(t,s)= \int\bds \eta_H(t-s,\nu) e^{2i\pi\nu t}\,d\nu
\end{equation*}
As a consequence, for $\kappa_H\in (S_0,\mathbf{L}^2, S_0')$, we also have 
$\eta_H\in (S_0,\mathbf{L}^2, S_0')$. In particular, this leads to the following expression for a weak evaluation of Gelfand triple operators: 
\begin{equation}\label{fo:OpEval1}
\langle K,L\rangle_{(\cBB, \cHH, \cBB')} = \langle \kappa (K ) ,\kappa (L) \rangle_{(S_0,\mathbf{L}^2,S_0' )} =  \langle \eta(K ) ,\eta (L) \rangle_{(S_0,\mathbf{L}^2,S_0' )}\end{equation}
For $L = g\otimes f^{\ast}$, i.e. the tensor product with kernel $\kappa (s,t) = g(s)\overline{f}(t)$ and spreading function $\eta (b,\nu ) = \cV_g f(b,\nu ) $, we thus have:
\begin{equation}\label{fo:OpEval2}
\langle K,g\otimes f^{\ast}\rangle_{(\cBB, \cHH, \cBB')} = \langle \kappa (K ) ,\kappa (g\otimes f^{\ast}) \rangle_{(S_0,\mathbf{L}^2,S_0' )} =  \langle \eta(K ) ,\cV_g f \rangle_{(S_0,\mathbf{L}^2,S_0' )}
\end{equation}
Let us also mention that the spreading function is related to the operator's Kohn-Nirenberg symbol via a \emph{symplectic Fourier transform}, which we define for later reference.
\begin{Def}\label{SympFT}
The symplectic Fourier transform is formally defined by
\[\cF_s F(t,\xi ) = \int_{\mathbb{P}}F(b,\nu )e^{-2\pi i (b \xi - t\nu )}dbd\nu\]\end{Def}
The symplectic Fourier transform is a self-inverse unitary automorphism of the Gelfand triple $(S_0,\mathbf{L}^2, S_0' )$.
We will make use of the following relation.
\begin{lemma}\label{Le:FSSTFT}Assume that $f_1,f_2,g_1,g_2\in \mathbf{L}^2 (\mathbb{R})$. Then
\[\cF_s (\cV_{g_1}f_1 \overline{\cV_{g_2} f_2})(x,\omega) = (\cV_{f_1}f_2 \overline{\cV_{g_2}g_1})(x,\omega ).\]\end{lemma}
\underline{Proof:} The analoguous statement  for the conventional (Cartesian) Fourier transform reads\ 
 $\cF (\cV_{g_1}f_1 \overline{\cV_{g_2} f_2})(x,\omega) = (\cV_{f_1}f_2 \overline{\cV_{g_2}g_1})(-\omega,x )$ and has been shown in~\cite[Lemma~2.3.2]{gr03-2}. The fact,  that 
\[ [\cF_s F](x,\omega) = \hat{F}(-\omega,x )\] completes the argument.\foorp\\
Recall that the spreading function of the product of operators corresponds to the twisted convolution of the operators' spreading function.
Assume  $K_1$ in $(\cBB, \cHH, \cBB') $ and $K_2$ in $(\cBB', \cHH, \cBB) $ then 

\begin{equation}\label{fo:opTwist}
	\eta (K_2\cdot K_1) = \eta(K_2)\natural \eta(K_1).
\end{equation}

The spreading function representation of operators provides an
interesting time-frequency implementation for operators, 
stated in the following proposition. It
turns out to be closely connected to the tools described in the
previous section, in particular twisted convolution and STFT.

\begin{proposition}
\label{prop:TwistRep}
Let $H$ be in $(\cBB, \cHH, \cBB') $, and let
$\eta = \eta_H$ be its spreading function in
$(S_0(\RR^2),\mathbf{L}^2(\RR^2), S_0'(\RR^2)$.
Let $g\in S_0(\RR)$,
then the STFT of $Hf$ is given by a twisted convolution of $\eta_H$ and $\cV_g f$:
\[\cV_g Hf (z) = (\eta_H\natural \cV_g f)(z).\]
\end{proposition}
\underline{Proof:} 
By~(\ref{fo:SFrep}), we may write
\bea
\cV_g Hf (z') &=&\langle Hf, \pi (z') g\rangle \notag\\
&=&\int \langle \eta_H(z)\pi (z) f, \pi (z' )g\rangle dz\label{SFrepeq1}\\
&=&\int \eta_H(z)\langle  f, \pi (z)^{\ast}\pi (z' )g\rangle dz\notag\\
&=&\int \eta_H(z)e^{-2\pi i b(\nu'-\nu )} \langle  f, \pi (z'-z )g\rangle dz\notag\\
&=&\int \eta_H(z)\cV_g f(z'-z) e^{-2\pi i b(\nu'-\nu )}  dz = (\eta_H\natural \cV_g f)(z')\notag
\eea
Note that $S_0$ is time-frequency shift-invariant, so $\pi (z) g$ is in $S_0$ for all $z$. Hence, the expression in \eqref{SFrepeq1} is well-defined. \\
If $f\in\mathbf{L}^2(\RR)$ and $H\in\cHH$, then $\cV_gf$,
$\cV_g Hf$ and $\eta_H \in\mathbf{L}^2(\RR^2)$, which is in
accordance with the fact that
$\mathbf{L}^2\natural \mathbf{L}^2\subseteq\mathbf{L}^2$. \\
If  $f\in S_0(\RR)$, then
$H$ may be in $\cBB'$, such that $\eta_H \in S_0'(\RR^2)$, hence $Hf\in S_0'(\RR )$. Hence, we have
$\cV_gf\in S_0(\RR^2)$ and $\cV_g Hf\in S_0'(\RR^2)$. This leads to the inclusion
$S_0'\natural \mathbf{W}(\mathbf{C},\ell^1)\subseteq \mathbf{L}^\infty$,
which may easily be verified directly. \\
On the other hand, if $f\in S_0'(\RR)$  and $H$  in $\cBB$, such
that $\eta_H \in S_0(\RR^2)$, then $Hf$ is in $S_0(\RR)$. Hence, 
$\cV_gf\in S_0'(\RR^2)$ and
$\cV_g Hf\in S_0(\RR^2)$. Here,  this leads to the
conclusion that we have,  for $f\in S_0'(\RR)$: 
\begin{equation}\label{eq:twistS0}
S_0\natural \cV_g f\subseteq S_0 .
\end{equation}
\foorp\\
Although it is known that $\cV_g f$ is not only
$\mathbf{L}^{\infty} (\RR^2 )$, but also in the Amalgam
space $\mathbf{W}(\mathcal{F}\mathbf{L}^1, \ell^{\infty})$ for
$f\in S_0'(\RR)$ and $g\in S_0(\RR)$,\cite{FZ98a}, it is
not clear, whether \eqref{eq:twistS0} also holds for
functions $F\in \mathbf{W}(\mathcal{F}\mathbf{L}^1, \ell^{\infty})$,
which are not in the range of $S_0' (\RR )$ under $\cV_g$.
This and other interesting open questions concerning the
twisted convolution of function spaces are currently under
investigation\footnote{H.~Feichtinger and F.~Luef.
Twisted convolution properties for {W}iener amalgam spaces.
{\em {I}n preparation}, 2008.}.

\begin{remark}\rm As a consequence of the last proposition, 
$H$ may be realized as a
twisted convolution in the time-frequency domain: 
\begin{equation}
\label{fo:TwistRep}
Hf = \int\bds\int\bds \left(\eta_H\natural \cV_gf \right)(b,\nu)
M_\nu T_b h\,dbd\nu\ \mbox{ for all }
f\in (S_0 , \mathbf{L}^2, S_0').
\end{equation}

Notice that Proposition~\ref{prop:TwistRep} implies that the range of $\cV_g$ is
invariant under left twisted convolution.
Notice also that this is no longer true if the
left twisted convolution is replaced with the right
twisted convolution. Indeed, in such a case, one has
$$
\cV_g f \natural \eta_H = \cV_{H^* g}f\ .
$$
Hence, one has the following simple rule: left twisted convolution on the
STFT amounts to acting on the analyzed function $f$, while right
twisted convolution on the STFT amounts to acting on the analysis window $g$.
It is worth noticing that in such a case, applying $\cV_g^*$ to
$\cV_gf \natural \eta_H$ yields the analyzed function $f$, up to some
(possibly vanishing) constant factor.
\end{remark}
\begin{example}As an illustrating example, let $g,h\in S_0$ be
such that $\langle g,h\rangle = 1$ and  consider the oblique projection $P: f\mapsto \langle f,g\rangle h$. The spreading function of this operator is given by $\cV_g h$, and we have $\cV_{\varphi} Pf (z)  = \langle f,g\rangle \langle h,\pi (z) \varphi\rangle$.  By virtue of the inversion formula for the STFT, which may be written as $\langle f,g\rangle h = \int \langle h , \pi (z) g\rangle \pi (z) f dz$, we obtain: 
\[\cV_{\varphi} Pf (z)  = \int\langle h,\pi (z) g\rangle\langle \pi (z) f, \pi (z') \varphi\rangle dz = \cV_g h\natural \cV_{\varphi}f.\]
By completely analogous reasoning, we obtain the converse formula, if the operator is applied to the analysing window: 
\[\cV_{P\varphi}f = \langle f, \pi (z) P\varphi\rangle = \cV_{\varphi}f\natural \cV_g h .\]
\end{example}
\begin{remark}\rm
Notice also that twisted convolution in the phase space is 
associated with the true translation structure. Indeed,
time-frequency shifts take the form of twisted convolutions
with a Dirac distribution on $\PP$:
$$
\delta_{b_0,\nu_0}\natural \cV_g f = \cV_g M_{\nu_0} T_{b_0} f\ .
$$
This corresponds to the usage of engineers, who  ``adjust the
phases'' after shifting STFT coefficients~\cite{dols1,phavoc66}.\end{remark}

\subsection{Time-Frequency  multipliers}\label{Se:tfmult}
\label{se:mult}
Section~\ref{Se:sprep}  has shown the close connection between the
spreading function representation of Hilbert-Schmidt operators and
the short time Fourier transform. However, the twisted convolution
representation is generally of poor practical interest in the
continuous case, because it does not discretize well. Even in the
finite case, it relies on the full STFT on $\CC^N$, which represents
vectors with $N^2$ STFT coefficients, which may be far too large
in practice, and sub-sampling is not possible in a
straightforward way.

Time-frequency (in particular Gabor) multipliers represent a valuable
alternative for time-frequency operator representation
(see~\cite{Feichtinger02first,Hlawatsch03linear} and references
therein for reviews). We analyze below the connections between these
representations and the spreading function, and point
out some limitations, before turning to generalizations.

\subsubsection{Definitions and main properties}
Let $g,h\in S_0(\RR)$ be such that $\langle g,h\rangle =1$,
let $\bm\in \mathbf{L}^\infty(\RR^2)$, and define the STFT multiplier
$\MM_{\bm;g,h}$ by 
\begin{equation}
\MM_{\bm;g,h} f =
\int_\PP \mathbf{m}(b,\nu) \cV_gf(b,\nu)\, \pi(b,\nu) h\, dbd\nu ,
\end{equation}
This defines a bounded operator on $( S_0(\RR), \mathbf{L}^2(\RR), S_0'(\RR))$.

Similarly,  given lattice constants $b_0,\nu_0\in\RR^+$, set
$\pi_{mn} = \pi(mb_0,n\nu_0) = M_{n\nu_0}T_{mb_0}$.
Then, for $\bm\in\ell^\infty(\ZZ^2)$, the
corresponding Gabor multiplier is defined as
\begin{equation}\label{fo:GMdef}
\MM_{\bm;g,h}^G f = \sum_{m=-\infty}^\infty
\sum_{n=-\infty}^\infty \bm(m,n)
\cV_gf(mb_0,n\nu_0)\,  \pi_{mn} h\ .
\end{equation}
Note that Gabor multipliers may be interpreted as STFT multipliers with multiplier $\bm$  in $S_0'$. In fact, in this case, $\bm$ is simply a sum of weighted Dirac impulses on the sampling lattice. 

The definition of time-frequency multipliers can of course
be given for $g,h\in \mathbf{L}^2 (\mathbb{R})$, many nice properties
only apply with additional assumptions on the windows.
Abstract properties of such multipliers have been studied extensively,
and we refer to~\cite{Feichtinger02first} for a review. One may show
for example that, whenever
the windows $g$ and $h$ are at least in $S_0$, if $\bm$ belongs to
$\mathbf{L}^2(\PP)$ (or $\ell^2(\ZZ^2)$) then the corresponding
multiplier is a Hilbert-Schmidt operator and maps $S_0'(\mathbb{R})$
to $\mathbf{L}^2 (\mathbb{R})$.

The spreading function of time-frequency multipliers may be computed
explicitly.
\begin{lemma}
The spreading function of the STFT multiplier $\MM_{\bm;g,h}$
is given by
\begin{equation}
\label{fo:STFTMult.spreading}
\eta_{\MM_{\bm;g,h}} (b,\nu) = \cM(b,\nu) \cV_g h(b,\nu)\ ,
\end{equation}
where $\cM$ is the symplectic Fourier transform of the transfer
function $\bm$
$$\cM(t,\xi) = \int_\PP \bm(b,\nu) e^{2i\pi(\nu t - \xi b)}\, dbd\nu\ .
$$
Specifying to the Gabor multiplier $\MM_{\bm;g,h}^G$, we see the same expression for the spreading function, however, in this case, 
$\cM = \cM^{(d)}$  is the $(\nu_0\inv,b_0\inv)$-periodic
symplectic Fourier transform of the discrete transfer function $\bm$
\begin{equation}
\label{fo:GabMult.spreading}
\cM^{(d)}(t,\xi) = \sum_{m=-\infty}^\infty\sum_{n=-\infty}^\infty
\bm (m,n) e^{2i\pi(n\nu_0 t - mb_0\xi)}\ .
\end{equation}
\end{lemma}
\underline{\em Proof~:}~For
$f\in S_0'$ and $\varphi\in S_0$, we may write
\[\langle \MM_{\bm;g,h} f,\varphi \rangle =
\langle \bm, \overline{\cV_g f}\cdot \cV_h\varphi\rangle\ ,
\]
where the right-hand side inner product has to be interpreted as an integral or infinite sum, respectively. 
By Lemma~\ref{Le:FSSTFT}, applying the symplectic Fourier transform, we obtain
\[
\langle \MM_{\bm;g,h} f,\varphi \rangle
= \langle \cM, \cV_{\varphi} f\cdot\overline{\cV_g h} \rangle
= \langle \cM \cdot\cV_g h, \cV_{\varphi} f\rangle\ .
\]
By calling on \eqref{fo:OpEval2}, this proofs~\eqref{fo:STFTMult.spreading}.
By virtue of that fact that for $g, h \in S_0$, $\cV_g h$ is certainly
in $\mathbf{L}^1 (\mathbb{P})$ and even in the Wiener Amalgam Space
$W(C, \mathbf{L}^1)$, hence in particular continuous, the expressions
for the spreading function given in the lemma are always
well-defined. The symplectic Fourier transform is a Gelfand triple isomorphism of $(S_0,\mathbf{L}^2, S_0')$, i.e.,  $\bm \in (S_0,\mathbf{L}^2, S_0') \Longleftrightarrow \cM\in (S_0,\mathbf{L}^2, S_0')$.  Hence,   $\cM\cdot\cV_g h\in (S_0,\mathbf{L}^2, S_0')$, which is in accordance with the fact, that for $\bm$ in $(\ell^1,\ell^2,\ell^{\infty})$, i.e., $(S_0,\mathbf{L}^2, S_0') (\ZZ^2)$, the kernel of  the resulting operator (and hence its spreading function), is in $(S_0,\mathbf{L}^2, S_0')$, see~\cite{Feichtinger02first} for details.
\foorp
\begin{remark}\rm
All  expressions derived so far are easily generalized to Gabor frames for $\mathbb{R}^d$ associated to
arbitrary lattices $\Lambda\subset\RR^{2d}$. In such situations, the spreading
function takes a similar form, and involves some discrete symplectic
Fourier transform of the transfer function $\bm$, which is in that case
a $\Lambda^\circ$-periodic function, $\Lambda^\circ$ being the adjoint
lattice of $\Lambda$, see Definition~\ref{Def:AdLatt}
in Section~\ref{se:multmult}.
\end{remark}

Notice that as a consequence of Theorem~\ref{th:SFrep},
 one has the following ``intertwining property''
$$
\cV_g \MM f = (\cM_{\bm;g,h}\,\cV_gh)\natural\cV_gf\ .
$$
\begin{remark}\rm
\label{rem:aaaaa}
It is clear from the above calculations that a general
Hilbert-Schmidt operator may not be well represented by a  
TF-multiplier.
For example, let us assume that the analysis and synthesis windows
have been chosen, and let $\eta$ be the spreading function
of the operator under consideration.
\begin{itemize}
\item In the STFT case, if the analysis and synthesis windows
are fixed, the decay of the spreading function has to be fast enough
(at least as fast as the decay of $\cV_gh$) to ensure the boundedness of
the quotient $\cM=\eta/\cV_gh$. {
Such considerations have led to
the introduction of the notion of {\em underspread}
operators~\cite{Kozek96matched} whose spreading function is
compactly supported in a domain of small enough area.} A more precise definition of underspread operators will be given below.
\item
In the Gabor case, the periodicity of $\cM^{(d)}$ imposes extra
constraints on the spreading function $\eta$. In particular, the shape
of the support of the spreading function must influence the choice of
optimal parameters for the approximation by a Gabor multiplier, i.e.
the shape of the window as well as the lattice parameters.
The following  numerical example indicates the direction for the
choice of parameters in the approximation of operators by
Gabor multipliers. 
\end{itemize}
\end{remark}
\begin{example}\label{Ex1}
Consider two operators $OP1$ and $OP2$ with spreading functions as shown
in Figure~\ref{FI1}. The values of the spreading functions are random
and real, uniformly distributed in $[-0.5,0.5]$.\\
Operator~1 has a spreading function with smaller support on the
time-axis, which means that the corresponding operator exhibits
time-shifts across smaller intervals than  Operator~2, whose
spreading function is, on the other hand, less extended in frequency.
The effect in the opposite direction is, obviously, reverse. These
characteristics are illustrated by applying the operators to a
sinusoid with frequency $1$ and a Dirac impulse at $-1$, respectively.\\ 
Next, we realize approximation\footnote{Best approximation is realized
in Hilbert Schmidt sense, see the next section for details.}  by
Gabor multipliers with two fixed pairs of lattice constants:
$b_1 = 2, \nu_1 = 8$ and $b_2 = 8, \nu_2 = 2$.  Furthermore, the windows are Gaussian windows
varying  from wide $(j = 0)$ to narrow $(j = 100)$.Thus, $j$ corresponds to the concentration of the window, in other words, $j$ is the reciprocal of the standard deviation.
Now the approximation
quality is investigated. The results are shown in the lower plots
of Figure~\ref{FI1}, where the left subplot shows the approximation
quality for operator $OP1$ for $b_1,\nu_1$ (solid) and
$b_2,\nu_2$ (dashed), while the right hand subplot gives the
corresponding results for operator~2. The error is measured by
$err = \|OP-APP\|_{\cHH}/\|OP+APP\|_{\cHH}$. Here, $APP$ denotes
the approximation operator and the norm is the operator norm. The
results show that, as expected, the "adapted" choice of time-frequency
parameters leads  to more favorable approximation quality.
Here, the adapted choice of $b$ and $\nu$ mimics the shape of
the support of the spreading function according to
formula~\eqref{fo:STFTMult.spreading} and the periodicity of
$\cM^{(d)}$. In brief, if the operator realizes frequency-shifts
in a wider range, we will need more sampling-points in frequency
and vice-versa. It is also visible, that the shape of the window
has considerable influence on the approximation quality.
\begin{figure}

\centerline{\includegraphics[scale = .95]{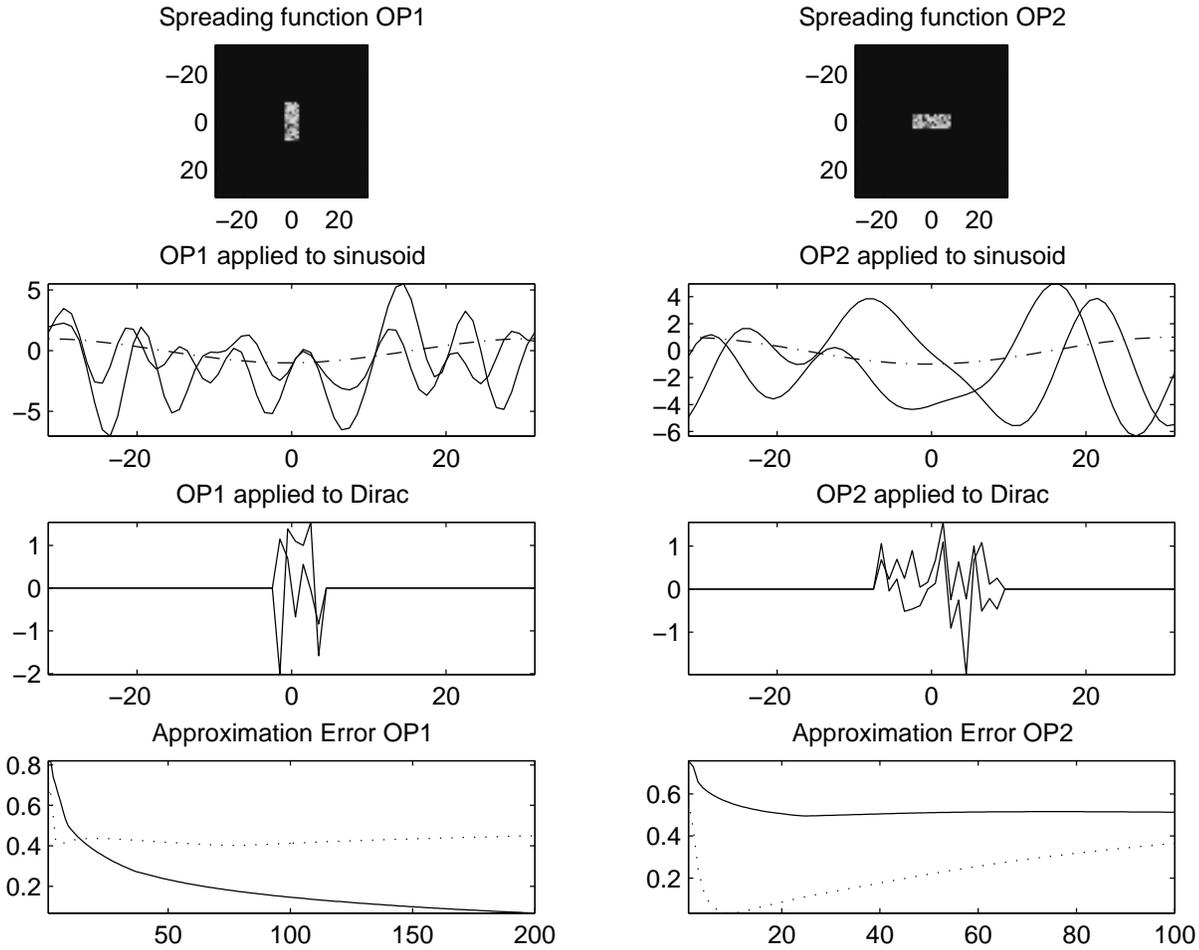}}
\caption{{\it Approximation by Gabor Multipliers with parameters $b_1, \nu_1$ (solid) and $b_2, \nu_2$ (dashed) }}\label{FI1}
\end{figure}
\end{example}
The previous example shows, that the parameters in the approximation by Gabor multipliers must be carefully chosen. Let us point out that the approximation quality achieved in the experiment described in Example~\ref{Ex1} is not satisfactory, especially when the time- and frequency shift parameters are not well adapted. Operators with a spreading function that is not well-concentrated around $0$, i.e. "overspread operators",   don't seem to be well-represented by a Gabor multiplier even with high redundancy (the redundancy used  in the example is $8$). Moreover, a realistic operator will have a spreading function with a much more complex shape. The next section will give some more details on approximation by Gabor multipliers before generalizations,  which allow for approximation of more complex operators,  are suggested.

\subsubsection{Approximation by Gabor multipliers}
The possibility of approximating operators by Gabor multipliers in
Hilbert-Schmidt sense depends on the properties of the rank one
operators associated with time-frequency shifted copies of the
analysis and synthesis windows.

Let $g,h\in S_0(\RR )$ be such that $\langle g,h\rangle = 1$. Let $\lambda =
(b_1,\nu_1)\in\PP$, and consider the rank one operator (oblique
projection) $P_\lambda$ defined by
\begin{equation}
P_\lambda f = (g_\lambda^{\ast}\otimes  h_\lambda) f  =
\langle f,g_\lambda\rangle h_\lambda\ ,\quad
f\in(S_0(\RR),\mathbf{L}^2(\RR),S_0'(\RR))\ .
\end{equation}
Direct calculations show that
the kernel of $P_\lambda$ is given by
\begin{equation}
\kappa_{P_\lambda} (t,s) = \overline{g}_\lambda(s) h_\lambda(t)\ ,
\end{equation}
and its spreading function reads
\begin{equation}
\eta_{P_\lambda} (b,\nu) = e^{2i\pi(\nu_1 b-b_1\nu)} \cV_gh(b,\nu)\ .
\end{equation}

The following result characterizes the situations for which
time-frequency rank one operators form a Riesz sequence, in which case
the best approximation by a Hilbert-Schmidt operator is well-defined. This result first appeared in~\cite{Feichtinger02wavelet}. Here, we give a slightly different version, which is obtained from the original statement by applying Poisson summation formula. This result was also given in~\cite{bepf06} for general full-rank lattices in $\mathbb{R}^d$.
\begin{proposition}\label{Prop:Ucond}
Let $g,h\in\mathbf{L}^2(\RR)$, with $\langle g,h\rangle\ne 0$,
let $b_0,\nu_0\in\RR^+$, and set
\begin{equation}
\cU(t,\xi) = \sum_{k,\ell=-\infty}^\infty
\left|\cV_gh\left(t+\frac{k}{\nu_0},\xi + \frac\ell{b_0}\right)\right|^2\ .
\end{equation}
The family
$\{P_{mb_0,n\nu_0},\,m,n\in\ZZ\}$ is a Riesz sequence in $\cH$ if and only if
there exist real constants $0<A\le B<\infty$ such that 
\begin{equation}
\label{fo:Ucond}
0 < A\le \cU(t,\xi)\le B<\infty\ \mbox{ a.e. on } [0,\nu_0\inv[\times [0,b_0\inv[ .
\end{equation}
We call this condition the $\cU$ condition.
\end{proposition}

It turns out, that the approximation of a given operator via a
standard minimization process yields an expression, which is only
well-defined if the  $\cU$ condition~(\ref{fo:Ucond}) holds.
\begin{theorem}
\label{th:gabmult.app}
Assume that  $\cV_gh$ and $b_0,\nu_0\in\RR^+$ are  such that
the $\cU$ condition~(\ref{fo:Ucond}) is fulfilled. Then the best
Gabor multiplier approximation (in Hilbert-Schmidt sense) of
$H\in\cHH$ is defined by the time-frequency transfer function $\bm$
whose discrete symplectic Fourier transform reads
\begin{equation}\label{Eq:MBestGM}
\cM(b,\nu) =
\frac{\sum_{k,\ell=-\infty}^\infty
\overline{\cV_gh}\left(b+k/\nu_0,\nu + \ell/b_0\right)
\eta_H\left(b+k/\nu_0,\nu + \ell/b_0\right)}
{\sum_{k,\ell=-\infty}^\infty
\left|\cV_gh\left(b+k/\nu_0,\nu + \ell/b_0\right)\right|^2}
\end{equation}
\end{theorem}
\underline{Proof: }
Let us denote as before by $\square$ the rectangle
$\square = [0,\nu_0\inv[\times [0,b_0\inv[$, and set $\cV=\cV_gh$
for simplicity of notation. First, notice that if $\eta_H\in L^2(\PP)$, then the function
$(b,\nu)\in\square\to \sum_{k,\ell} |\eta_H(b+k/\nu_0,\nu+\ell/b_0)|^2$
is in $L^2(\square)$, and is therefore well defined almost everywhere
in $\square$. Thus, by Cauchy-Schwarz inequality, the numerator
in~\eqref{Eq:MBestGM} is well-defined a.e.\\
The Hilbert-Schmidt optimization is equivalent to the problem
$$
\min_{\cM\in {\mathbf L}^2(\square)} \|\eta_H -\cM\cV\|^2\ .
$$
The latter squared norm may be written as
\begin{eqnarray*}
\|\eta_H\! -\!\cM\cV\|^2 &=&\!\!
\int\bds\int\bds \left|\eta_H(b,\nu) -\cM(b,\nu)\cV(b,\nu)\right|^2\,dbd\nu\\
&=&\!\!\!\! \sum_{k,\ell=-\infty}^\infty\!\!\int\!\!\!\int_\square
\left|\eta_H(b+k/\nu_0,\nu+\ell/b_0)-\cM(b,\nu)\cV(b+k/\nu_0,\nu+\ell/b_0)\right|^2
\,dbd\nu\\
&=&\!\! \int\!\!\!\int_\square \bigg[\sum_{k,\ell}
|\eta_H(b+k/\nu_0,\nu+\ell/b_0)|^2\\
&&\hphantom{aa}- 2 \Re\bigg( \overline{\cM}(b,\nu) \sum_{k,\ell}
\eta_H(b+k/\nu_0,\nu+\ell/b_0)\overline{\cV}(b+k/\nu_0,\nu+\ell/b_0)\bigg)\\
&&\hphantom{aa} +
|\cM(b,\nu)|^2 \sum_{k,\ell} |{\cV}(b+k/\nu_0,\nu+\ell/b_0)|^2
\bigg]db d\nu
\end{eqnarray*}
From this expression, the Euler-Lagrange equations may be obtained, which read
$$
\cM(b,\nu) \sum_{k,\ell} |{\cV}(b+k/\nu_0,\nu+\ell/b_0)|^2
= \sum_{k,\ell}
\eta_H(b+k/\nu_0,\nu+\ell/b_0)\overline{\cV}(b+k/\nu_0,\nu+\ell/b_0)\ ,
$$
and the result follows.\foorp\\

We  next derive an error estimate for the approximation. Let us set, for $(b,\nu)\in\square$,
$$
\cE(b,\nu) = \frac{\left|\sum_{k,\ell}\eta_H(b+k/\nu_0,\nu+\ell/b_0)
\overline{\cV}(b+k/\nu_0,\nu+\ell/b_0)\right|^2}
{\sum_{k,\ell}|\eta_H(b+k/\nu_0,\nu+\ell/b_0)|^2\,\cU(b,\nu)}\  .
$$
\begin{corollary}\label{Cor:ErrEst1}
With the above notation, we obtain the estimate
\begin{equation*}
\|H-\MM_\bm\|^2_\cHH
\le \|\eta_H\|^2\,\|1-\cE\|_\infty
\end{equation*}for the best approximation of $H$ by a Gabor multiplier $\MM_\bm$ according to \eqref{Eq:MBestGM}.
\end{corollary}

\underline{\em Proof~:}
Set $\Gamma_H(b,\nu) = \sum_{k,\ell}|\eta_H(b+k/\nu_0,\nu+\ell/b_0)|^2$. 
Replacing the expression for $\cM$ obtained in \eqref{Eq:MBestGM}
into the error term, we have
\begin{eqnarray*}
\|H-\MM_\bm\|^2_\cHH
&=& \int\int_\square \left[\Gamma_H(b,\nu)
 - |\cM(b,\nu)|^2\cU(b,\nu)\right]\,dbd\nu\\
&=&\int\int_\square \Gamma_H(b,\nu)\left[1-\cE(b,\nu)\right]\,dbd\nu\\
&\le& \|\eta_H\|^2\,\|1-\cE\|_\infty
\ ,
\end{eqnarray*}
where we have used the fact that $\|\Gamma_H\|_{{\mathbf L}^2(\square)}
= \|\eta_H\|_{{\mathbf L}^2(\PP)}$. 
Clearly, Cauchy-Schwarz inequality gives $|\cE(b,\nu)|\le 1$
on $\square$, with equality
if and only if there exists a function $\phi$ such that
$\eta = \phi\cV$, i.e. if and only if $H$ is a multiplier with
the prescribed window functions. Hence we  obtain
\begin{equation}
\label{fo:GM.error}
\|H-\MM_\bm\|^2_\cHH \le \|H\|_\cHH^2\,
\left[1-\mathop{{\mathrm ess}\inf}\limits_{(b,\nu)\in\square} \cE(b,\nu)\right]
\ .
\end{equation}\foorp
\begin{remark}\rm

$\cE(b,\nu)$ essentially represents the cosine of the angle between
vectors $\{\eta_H(b+k/\nu_0,\nu+\ell/b_0),\,k,\ell\in\ZZ\}$
and $\{\cV(b+k/\nu_0,\nu+\ell/b_0),\,k,\ell\in\ZZ\}$. In other words, the
closer to colinear these vectors, the better the approximation.\end{remark}
}
An
example for operators which are poorly represented by this class of
multipliers are those with a spreading function that is not
``well-concentrated''. These are,  in technical terms, overspread
operators.
The underspread/overspread terminology seems to originate from the context of time-varying multipath wave propagation channels~\cite{Ken69}. However, different definitions exist in the literature. Here, we give the definition used in~\cite{ko97-1}. Note that underspread operators have recently found renewed interest~\cite{ko97,bokoma02,hlma98,kopf06,pfwa05-1}.  
\begin{Def}
Consider  an operator $K$ with compactly supported spreading function:
\[\supp(\eta_K )\subseteq Q(t_0,\xi_0), \mbox{ where }
Q(t_0,\xi_0):= [-t_0,t_0]\times [-\xi_0,\xi_0].\]
Then, $K$ is called underspread, if $t_0\xi_0<1/4$.
\end{Def}
Most generally, operators which are \emph{not} underspread, will be
called overspread.\\
It is generally known, that an operator must be underspread in order to be well-approximated by  Gabor multipliers. Formula \eqref{Eq:MBestGM} enables us to make this statement more precise.
\begin{corollary}
Consider an underspread operator $H$. Then, for $(b_0,\nu_0)$ such that
$1/b_0>2\xi_0$ and $1/\nu_0>2t_0$, it is possible to find a
Gabor frame $\{g_{mn},m,n\in\ZZ\}$, with lattice constants
$(b_0,\nu_0)$, and dual window $h$.
Then,  the symplectic Fourier
transform $\cM$ of the time-frequency transfer function of the
best Gabor multiplier takes the  form
$$
\cM = \frac{\eta_H\overline{\cV_gh}}{\cU}\ ,
$$
and the approximation error can be bounded by 
\begin{equation*}
\left\|\eta_H-\cM\cV_gh\right\|^2 \le \|\eta_H\|^2
\mathop{{\mathrm ess}\sup}\limits_{(b,\nu)\in\square}
\left[ 1 - \frac{|\cV_gh(b,\nu)|^2}{\cU(b,\nu)}\right]
\end{equation*}
\end{corollary}
\underline{Proof: }
\begin{eqnarray*}
\left\|\eta_H-\cM\cV_gh\right\|^2 &=&
\int\int_\square \left[|\eta(b,\nu)|^2 
-\frac{\left|\eta(b,\nu)\cV_gh(b,\nu)\right|^2}{\cU(b,\nu)}\right]\,dbd\nu\\
&\le& \|\eta_H\|^2 \mathop{{\mathrm ess}\sup}\limits_{(b,\nu)\in\square}
\left[ 1 - \frac{|\cV_gh(b,\nu)|^2}{\cU(b,\nu)}\right]
\end{eqnarray*}
\foorp\\
 The  estimate in the last corollary shows, that approximation quality is  a joint property of window and lattice, which is in accordance with the results of Example~\ref{Ex1}. 
\begin{remark}\rm
Note that, although technically only defined for Hilbert-Schmidt
operators, the approximation by Gabor multipliers can formally be
extended to operators from $\cBB'$,
see~\cite[Section~5.8]{Feichtinger02first}. Also, the expression given in \eqref{Eq:MBestGM} is well-defined in $S_0'$ whenever $\eta_H$ is at least in $S_0'$. However, for non-Hilbert-Schmidt operators, it is not clear, in which sense the resulting Gabor multiplier represents the original operator. The following example shows that at least in some cases, the result is however the intuitively expected one.
\end{remark}
\begin{example}
Consider the operator $\pi (\lambda )$, i.e. a time-frequency shift.
Although this operator is clearly not a Hilbert-Schmidt operator,
we may consider its   approximation by a Gabor-multiplier according
to~\eqref{Eq:MBestGM}. First note that the spreading function of
the time-frequency shift $\pi (\lambda )= \pi (b_1,\nu_1)$ is given
by $\eta_{\pi} = \delta (b-b_1)\cdot \delta (\nu -\nu_1)$.
Then, we have
\begin{eqnarray*}	
\cM (b,\nu ) &= &
\sum_{k,l}\frac{\overline{\cV_g h}(b+k/\nu_0, \nu+l/b_0 )\delta(b-b_1+k/\nu_0, \nu-\nu_1+l/b_0)}{\sum_{k',l'}|\cV_g h(b+k'/\nu_0, \nu+l'/b_0 )|^2}\\
&=&
\overline{\cV_g h }(b_1,\nu_1)\sum_{k,l}\frac{\delta(b-b_1+k/\nu_0, \nu-\nu_1+l/b_0)}{\sum_{k',l'}
|\cV_g h(b_1+\frac{k'-k}{\nu_0},\nu_1+\frac{l'-l}{b_0})|^2}\\
&=&
\frac{\overline{\cV_g h }(b_1,\nu_1)}{\cU(b_1,\nu_1)}
\sum_{k,l}\delta(b-b_1+k/\nu_0, \nu-\nu_1+l/b_0)
\end{eqnarray*}	

Hence, from the inverse (discrete) symplectic Fourier transform we obtain:
\begin{eqnarray*}
\bm (m,n)& = &
\frac{\overline{\cV_g h }(b_1,\nu_1)}{\cU(b_1,\nu_1)}
\int_0^{\frac{1}{\nu_0}}\int_0^{\frac{1}{b_0}}
\sum_{k,l}\delta(b-b_1+k/\nu_0, \nu-\nu_1+l/b_0)\\
&&
\qquad\qquad\qquad\qquad\qquad\times
 e^{-2\pi i (mb_0 \nu  -n\nu_0 b)}\,dbd\nu\\
&=&
\frac{\overline{\cV_g h }(b_1,\nu_1)}{\cU(b_1,\nu_1)}\ 
e^{-2\pi i (mb_0 \nu_1  -n\nu_0 b_1)}
\end{eqnarray*}	
As expected, the absolute value of the mask is constant and the
phase depends on the displacement of $\lambda$ from the origin.
{
 This confirms the key role played by the phase
of the mask of a Gabor multiplier.}
Specializing to $\lambda = 0$, we obtain a constant mask and
thus,  if $h$ is a dual window of $g$ with respect to $\Lambda = b_0\mathbb{Z}\times \nu_0\mathbb{Z}$, up to a constant
factor, the identity.
\end{example}

\section{Generalizations: multiple Gabor multipliers and TST spreading
functions}
\label{se:multmult}
In the last section it has become clear that most operators are not
well represented as a STFT or Gabor multiplier. 

Guided  by the desire to extend the good approximation quality that
Gabor multipliers warrant for underspread operators to the
class of their overspread counterparts, we introduce
generalized TF-multipliers. The basic idea is to allow for an
extended scheme in the synthesis part of the operator: instead of
using just one window $h$, we suggest the use  of a set of windows
$\{h^{(j)}\}$ in order to obtain the class of 
{\em Multiple Gabor Multipliers} (MGM for short).

\begin{Def}[Multiple Gabor Multipliers]
Let  $g\in S_0 (\mathbb{R})$ and a family of reconstruction
windows $h^{(j)}\in S_0 (\mathbb{R})$, $j\in \mathcal{J}$, as well
as corresponding masks $\bm_j\in\ell^{\infty}$ be given.
Operators of the form
\begin{equation}
\label{fo:multimulti}
\MM = \sum_{j\in\cJ} \MM^G_{\bm_j;g,h^{(j)}}
\end{equation}
will be called  {\em Multiple Gabor Multipliers} (MGM for short).
\end{Def}
Note, that  we need to impose additional assumptions in order
to obtain a well-defined operator.  For example, we may  assume
$\sum_j sup_{\lambda} |m_j (\lambda )| = C< \infty$ and
$\max_j \|h^j\|_{S_0} = C <\infty$, which guarantees a  bounded
operator on $(S_0,\mathbf{L}^2,S_0')$. This follows easily from
the boundedness of a Gabor multiplier under the condition that
$\bm$ is $\ell^{\infty}$. Conditions for function space membership
of MGMs are easily derived in analogy to the Gabor multiplier case.
For example, if $\sum_j\sum_{\lambda} |m_j(\lambda )|^2<\infty$, we
obtain a Hilbert-Schmidt operator, similarly, trace-class membership
follows from an analogous $\ell^1$-condition.

As a starting point, we give the (trivial) generalization of the
spreading function of a MGM as a sum of the spreading functions
corresponding to the single Gabor multipliers involved. 
\begin{lemma}
The spreading function of a MGM is (formally) given by
\begin{equation}
\label{fo:MGM.spreading}
\eta_{\MM_{\bm_j;g,h^{(j)}}}^G (b,\nu) =
\sum_{j\in\cJ}\cM^{(j)}(b,\nu) \cV_g h^{(j)}(b,\nu)\ ,
\end{equation}
where the $(\nu_0\inv,b_0\inv)$-periodic functions $\cM^{(j)}$ are the
symplectic Fourier transforms of the transfer functions $\bm_j$.
\end{lemma}
Note that the issue of convergence for the series defining
$\eta_{\MM_{\bm_j;g,h^{(j)}}}^G (b,\nu)$ will not be discussed, as
in practice $|\mathcal{J}|$ will usually be finite. Let us just
mention that by assuming $\bm_j\in\ell^2(\mathbb{Z}^3 )$, i.e.,
the Hilbert-Schmidt case with an additional $\ell^2$-condition
for the masks in the general model, we have
$\|\eta_{\MM_{\bm_j;g,h^{(j)}}}^G \|_2 = C\|\bm\|_2$.\\
It is immediately obvious that this new model gives much more
freedom in generating overspread operators. However, in order to
obtain structural results, we will have to impose further specifications.\\
Before doing so, we will state a generalization of
Proposition~\ref{Prop:Ucond} to the more general situation of the
family of projection operators  $P^j_\lambda$ defined by 
\begin{equation}
P^j_\lambda f = (g_\lambda^{\ast}\otimes  h^j_\lambda) f =
\langle f,g_\lambda\rangle h^j_\lambda , \mbox{ where }
\lambda\in\Lambda, j\in\mathcal{J}, \ |\mathcal{J}|<\infty .
\end{equation}
Note that these  projection operators are the building blocks for
the MGM. The following theorem characterizes their Riesz property.

\begin{proposition}
\label{Prop:GenUcond}
Let $g,h^j\in\mathbf{L}^2(\RR)$,
$j\in \mathcal{J}, \ |\mathcal{J}|<\infty$, with
$\langle g,h^j\rangle\neq 0$,
let $b_0,\nu_0\in\RR^+$, and let the matrix
$\Gamma(b,\nu)$ be defined by
\begin{equation}
\label{EqU1}
\Gamma(b,\nu)_{jj'} =
\sum_{k,\ell}\overline{\cV_gh^{(j)}}(b+k/\nu_0,\nu+\ell/b_0)
\cV_gh^{(j')}(b+k/\nu_0,\nu+\ell/b_0)
\end{equation}
a.e. on $\square  = [0,\nu_0^{-1}[\times [0,b_0^{-1}[$.
Then the family of projection operators
$\{P^j_\lambda, j\in \mathbb{Z}, \lambda\in\Lambda\}$ is
a Riesz sequence  in $\cH$ if and only if $\Gamma$ is invertible a.e.\\
Alternatively, the Riesz basis property is characterized by
invertibility of the matrix $U$ defined as
\begin{equation}
\label{EqU2}
\cU^{jj'}(t,\xi)
= \sum_{k,l} U^{jj'}(kb_0,l\nu_0 ) e^{-2\pi i (l\nu_0 t-kb_0\xi)}
\end{equation}
a.e. on the fundamental domain of $\Lambda$.
\end{proposition} 
\underline{Proof:} Recall  that the family
$\{P^j_\lambda, j\in \mathbb{Z}, \lambda\in\Lambda\}$ is a Riesz
basis for its closed linear span 
if there exist constants $0<A,B<\infty$ such that 
\begin{equation}
A\|c\|_2^2\leq
\left\|\sum_{\lambda}\sum_j c^j_{\lambda} P^j_{\lambda}\right\|_{\HS}^2\leq
B\|c\|_2^2
\end{equation}
for all finite sequences $c$ defined on $(\Lambda\times \mathcal{J})$. We have
\begin{eqnarray*}
\left\|\sum_{\lambda}\sum_j c^j_{\lambda} P^j_{\lambda}\right\|_{\HS}^2
&=& \left\langle \sum_{\lambda}\sum_j  c^j_{\lambda} P^j_{\lambda}  ,
\sum_{\mu}\sum_{j'}c^{j'}_{\mu} P^{j'}_{\mu}\right\rangle\\
&=&
\sum_{\lambda\mu}\sum_{jj'} c^j_{\lambda}\overline{c^{j'}_{\mu}}
\overline{\langle g_{\lambda},g_{\mu} \rangle}
\langle h^j_{\lambda},h_{\mu}^{j'}  \rangle.\\
\end{eqnarray*}
Hence, by setting $U^{jj'} (b,\nu ) = [\overline{\cV_g g}\cdot \cV_{h^j}h^{j'}](b,\nu )$, we may write
\begin{eqnarray*}
\left\|\sum_{\lambda}\sum_j c^j_{\lambda} P^j_{\lambda}\right\|_{\HS}^2
&=&
\sum_{\lambda\mu}\sum_{jj'} c^j_{\lambda}\overline{c^{j'}_{\mu}}
U^{jj'} (\mu - \lambda )\\
&=&
\sum_{\mu}\sum_{jj'}\overline{c^{j'}_{\mu}}(c^j\ast U^{jj'})(\mu)\\
&=&
b_0\nu_0\sum_{jj'}\left\langle \cU^{jj'}\cdot\cC^j,
\cC^{j'}\right\rangle_{\mathbf{L^2}(\square )},
\end{eqnarray*}
where
$\cU^{jj'}(t,\xi) = \sum_{k,l} U^{jj'}(kb_0,l\nu_0 )
e^{-2\pi i (l\nu_0 t-kb_0\xi)}$ is the discrete symplectic Fourier
transform of $U^{jj'}$, and, analogously, $\cC^j$ is the discrete
symplectic Fourier transform of the sequence $c^j$, defined on
$\Lambda$, for each $j$. Hence, these are
$\nu_0^{-1}\times b_0^{-1}$-periodic functions.
The last equation can be rewritten as
 \begin{eqnarray*}
\left\|\sum_{\lambda}\sum_j c^j_{\lambda} P^j_{\lambda}\right\|_2^2
&=&
b_0\nu_0\sum_{jj'}\int_{\square} \cU^{jj'}(t,\xi )\cdot\cC^j(t,\xi )
\overline{\cC^{j'}}(t,\xi )dtd\xi\\
&=&
\left\langle \cU\cdot \cC,\cC\right\rangle_{\mathbf{L^2}(\square )
\times\ell^2(\mathbb{Z})},
\end{eqnarray*}
where $\cU$ is the  matrix with entries $\cU^{jj'}$.
Note that this proves statement \eqref{EqU2} by positivity of
the operator $\cU$.\\
In order to obtain the condition for $\Gamma$ given in~\eqref{EqU1},
first note that 
\[
\cF_s ( \overline{\cV_g g}\cdot \cV_{h^j}h^{j'})(\lambda)
= (\cV_g h^{j'}\cdot  \overline{\cV_g h^j})(\lambda)
\]
by applying Lemma~\ref{Le:FSSTFT}. Furthermore,
$F:=\cV_g h^{j'}\cdot  \overline{\cV_g h^j}$ is always in
$\mathbf{L}^1$ for $g,h^j\in \mathbf{L}^2$. We may therefore look
at the Fourier coefficients of its $\Lambda^{\circ}$-periodization,
with $\lambda = (mb_0 ,n\nu_0 )$:
\begin{eqnarray*}
\cF_s^{-1} (P_{\Lambda^{\circ}}F) (\lambda)
&=& \int_{\square}( \sum_{k,\ell}F(b+\frac{k}{\nu_0},\nu+\frac{\ell}{b_0}))
e^{2\pi i (b n\nu_0-\nu m b_0)}dbd\nu\\
&=&\int_{\square}( \sum_{k,\ell}F(b+\frac{k}{\nu_0},\nu+\frac{\ell}{b_0}))
e^{2\pi i ((b+\frac{k}{\nu_0})n\nu_0-(\nu+\frac{\ell}{b_0})m b_0 )}dbd\nu\\
&=&\int_{\mathbb{R}^2}F(b,\nu) e^{2\pi i (bn\nu_0-\nu mb_0)}dbd\nu
= \cF_s^{-1} (\cV_g h^{j'}\cdot  \overline{\cV_g h^j})(\lambda).
\end{eqnarray*}
Hence, we may apply the Poisson summation formula, with
convergence in $\mathbf{L^2}(\square )$, to obtain:
\begin{eqnarray*}
P_{\Lambda^{\circ}}(\cV_g h^{j'}\cdot  \overline{\cV_g h^j}) (b,\nu)
&=&
b_0\nu_0\sum_{k,l}\cF_s^{-1} (\cV_g h^{j'}\cdot 
\overline{\cV_g h^j})(kb_0,l\nu_0)e^{-2\pi i (b\ell\nu_0 -\nu kb_0)}\\
&=&
b_0\nu_0\sum_{k,l}   \overline{\cV_g g}\cdot \cV_{h^j}h^{j'}
(kb_0,l\nu_0)e^{-2\pi i (b\ell\nu_0 -\nu kb_0)}.
\end{eqnarray*}We conclude that 
\begin{equation}
\left\|\sum_j\sum_{\lambda\in\Lambda} c_{\lambda}^j P_{\lambda}^j\right\|_2^2
= b_0\nu_0\int\Gamma(b,\nu)_{jj'}C^j(b,\nu ) \overline{C^{j'} (b,\nu )}dbd\nu,
\end{equation}
and the Riesz basis property is equivalent to the invertibility of $\Gamma$.
\hfill\foorp\\

In the sequel, the discrete symplectic Fourier
transforms of $\bm_j$ will be denoted by $\cM_j$,
and the vector with $\cM_j$ as coordinates will be denoted by $\cM$. We then obtain an expression for the best multiplier in analogy to the Gabor multiplier case discussed in Theorem~\ref{th:gabmult.app}. 

\begin{proposition}
\label{prop:multgabmult.app}
Let $g\in S_0(\RR)$ and $h^{(j)}\in S_0(\RR)$, $j\in\mathcal{J}$ be such that
for almost all $b,\nu$, the matrix $\Gamma(b,\nu)$ defined in \eqref{EqU1} 
is invertible a.e. on $\square  = [0,\nu_0^{-1}[\times [0,b_0^{-1}[$.

Let $H\in (\cBB,\cHH,\cBB')$ be an operator with spreading function
$\eta\in (S_0, \mathbf{L}^2, S_0')$.
Then the functions $\cM_j$ yielding  approximation of
the form~(\ref{fo:multimulti}) may be obtained as 
\begin{equation}
\label{fo:tototo0}
\cM = \Gamma^{-1}\cdot \mathcal{B}\ ,
\end{equation}
where $\mathcal B$ is the vector whose entries read
\begin{equation}
\cB_{j_0}(b,\nu) 
= \sum_{k,\ell}\eta(b+k/\nu_0,\nu+\ell/b_0)
\overline{\cV}_g h^{j_0}(b+k/\nu_0,\nu+\ell/b_0).
\label{fo:tototo}
\end{equation}
For operators in $\cHH$ the obtained approximation is optimal in Hilbert-Schmidt sense.
\end{proposition}
\underline{Proof:}
The proof follows the lines of the Gabor multiplier case.
The optimal approximation of the form~(\ref{fo:multimulti}), when it
exists, is obtained by minimizing
$$
\big\|\eta - \sum_j \cM_j\cV_j\big\|^2
= \sum_{k,\ell}\int_\square \big|\eta(b+k/\nu_0,\nu+\ell/b_0)
- \sum_j \cM_j(b,\nu)\cV_j(b+k/\nu_0,\nu+\ell/b_0)\big|^2\,dbd\nu
$$
where one has set $\cV_j = \cV_gh^{(j)}$. Setting to zero the G\^ateaux
derivative with respect to $\overline{\cM}_{j_0}$, we obtain the
corresponding  variational equation
$$
\sum_j\cM_j(b,\nu)
\sum_{k,\ell}\cV_j(b+k/\nu_0,\nu+\ell/b_0)
\overline{\cV}_{j_0}(b+k/\nu_0,\nu+\ell/b_0) = \cB_{j_0}(b,\nu)\ ,
$$
where $\cB_j(b,\nu)$ are as defined in~(\ref{fo:tototo}).
Provided that the $\Gamma(b,\nu)$ matrices are invertible for almost all
$b,\nu$, this implies that the functions $\cM_j$ for 
approximation of the form~(\ref{fo:multimulti}) may indeed be obtained
as in~(\ref{fo:tototo0}).
\hfill\foorp.

In a next step, we are  going to  discern two basic approaches: \\
(a) $m_j (\lambda) = m(\mu,\lambda)$, i.e.~ the synthesis windows are time-frequency shifted versions (on a lattice) of a single synthesis window:
{
$h_j = \pi (\mu_j ) h$, $\mu_j\in\Lambda_1$.}\\
(b) $m_j (\lambda ) = m_1(\lambda)m_2(j)$, i.e.~ a separable multiplier function. 
If we set $h^{(j)} (t) = \pi(b_j,\nu_j)h(t)$ then this approach leads to what will be called TST spreading functions in Section~\ref{Se:TSTSF}. \\

In both cases we will be especially  interested in the situation in which the $h^j$ are given as time-frequency shifted versions of a single synthesis window on the  adjoint lattice $\Lambda^{\circ}$.
\begin{Def}[Adjoint lattice]\label{Def:AdLatt}
For a given lattice $\Lambda= b_0\mathbb{Z}\times \nu_0\mathbb{Z}$ the adjoint lattice  is given by $\Lambda^{\circ}=\frac{1}{\nu_0}\mathbb{Z}\times \frac{1}{b_0}\mathbb{Z}$.\end{Def} Note that the adjoint lattice is the dual lattice $\Lambda^{\perp}$ with respect to the symplectic character.\\
%
\subsection{Varying the multiplier: MGM with synthesis windows on the lattice}

We fix the synthesis windows  $h_j$ to be 
time-frequency translates of a fixed window function, i.e.
\begin{equation}
\label{fo:glou}
h^{(j)} (t) = \pi(b_j,\nu_j)h(t) = e^{2i\pi\nu_j t}h(t-b_j)\ .
\end{equation}
We may turn our attention to the projection operators associated to the (Gabor) families $(g,\Lambda_1)$ and $(h,\Lambda_2)$.
Note that it has been shown by Benedetto and Pfander~\cite{bepf06} that the family of projection operators $\{P_{\lambda},\,\lambda\in\Lambda\}$, as discussed in Section~\ref{se:mult} either 
forms a Riesz basis or not a frame (for its closed linear span). The next corollary shows that, on the other hand, if we use the extended family of projection operators $\{P_{\lambda,\mu}$, $(\lambda,\mu) \in\Lambda_1\times\Lambda_2\}$, where $P_{\lambda,\mu}f = \langle f,\pi (\lambda )g\rangle\pi (\mu ) h$, we obtain a frame of operators for the space of Hilbert-Schmidt operators, whenever $(g,\Lambda_1)$ and $(h,\Lambda_2)$ are Gabor frames. This corollary is a special case of Theorem~4.1 in~\cite{ba08a} and Proposition~3.2 in~\cite{ba08b}.
\begin{corollary}\label{Co:ProdFram}
Let two Gabor frames $(g,\Lambda_1)$ and $(h,\Lambda_2)$ be given. Then the family of projection operators $\{P_{\lambda,\mu}$, $(\lambda,\mu) \in\Lambda_1\times\Lambda_2\}$ form a frame of operators in $\cH$ and any Hilbert-Schmidt  operator $H$ may be expanded as
\[H  = \sum_{\lambda\in\Lambda_1, \mu\in\Lambda_2}\mathbf{c}(\lambda,\mu) P_{\lambda,\mu} .\]
The coefficients are given by $\mathbf{c}(\lambda,\mu) = \langle H, (P_{\lambda,\mu})^{\ast}\rangle =  \langle H \pi(\mu)h,\pi(\lambda)g\rangle$.
\end{corollary}
Note that an analogous statement holds for Riesz sequences. Very recently, it has been shown~\cite{bo08}, that the converse of Corollary~\ref{Co:ProdFram} holds true for both frames and Riesz bases, i.e. the family of projection operators $\{P_{\lambda,\mu}$, $(\lambda,\mu) \in\Lambda_1\times\Lambda_2\}$ is a frame (a Riesz basis) for  $\cH $ if and only if the two generating sequences form a frame (a Riesz basis) for $\mathbf{L}^2 (\RR)$. In particular, this leads to the conclusion, that the characterization of Riesz sequences given in Proposition~\ref{Prop:GenUcond} also  yields a characterization of frames for $\mathbf{L}^2 (\RR)$ -
it is well known, that $(g_{\lambda}, \lambda\in\Lambda)$ form a Gabor frame if and only if $(g_{\mu}, \mu\in\Lambda^{\circ})$ form a Riesz sequence. We can draw two conclusions.
\begin{corollary}\label{Co:ProdFram2}Let $g\in S_0 (\RR )$ and a lattice $\Lambda= b_0 \mathbb{Z}\times\nu_0\mathbb{Z}$ be given. 
\begin{itemize}\item[(a)] The Gabor family $\{g_{\lambda}, \lambda\in\Lambda\} $ forms a frame for $\mathbf{L}^2(\RR )$ if and only if the matrix
\begin{eqnarray*}\Gamma_{mn,m'n'}(b,\nu) 
&=& \sum_{k,\ell}
\exp\big(-2i\pi [m/\nu_0(\nu+\ell\nu_0 -n/b_0)- m'/\nu_0(\nu+\ell\nu_0 -n'/b_0)]\big)
\\
&&\hphantom{aaaaaa}\times\overline{\cV_g g}(b-m/\nu_0+kb_0,\nu-n/b_0+\ell\nu_0)\\
&&\hphantom{aaaaaa}\times{\cV_g g}(b-m'/\nu_0+kb_0,\nu-n'/b_0+\ell\nu_0)
\end{eqnarray*}is, a.e. on $\square$,  invertible on $\ell^2$.\\
\item[(b)] In addition, we may state the following "Balian-Low Theorem for the tensor products of Gabor frames":\\
A family of projection operators  given by $\{P_{\lambda,\mu} = g_\lambda^{\ast}\otimes  g_\mu$, $(\lambda,\mu) \in\Lambda\times\Lambda^{\circ}\}$ forms a frame for the space of Hilbert-Schmidt operators on $\mathbf{L}^2(\RR )$ if and only if it forms a Riesz basis. Hence, in this case, $g$ cannot be in $S_0 (\RR )$. 
\end{itemize}
\end{corollary}
\underline{Proof:} Statement (a) is easily obtained from \eqref{EqU1} by observing that $$
\cV_g \pi(mb_0,n\nu_0 )g(b,\nu) = e^{-2i\pi mb_0,(\nu-n\nu_0)}\cV_gh(b-mb_0,\nu-n\nu_0)\ .$$
We then have that $\{P_{\lambda,\mu} = g_\lambda^{\ast}\otimes  g_\mu$, $(\lambda,\mu) \in\Lambda^{\circ}\times\Lambda^{\circ}\}$ forms a Riesz basis in  $\cH$ if and only if 
$\Gamma_{mn,m'n'}(b,\nu) $ is invertible. By the converse of Corollary~\ref{Co:ProdFram}, this is equivalent to the Riesz property of $g_{\mu}, \mu\in\Lambda^{\circ} $, which, in turn, is equivalent to the frame property of $\{g_{\lambda}, \lambda\in\Lambda\} $ by Ron-Shen duality, see, e.g.~\cite{gr01}.\\
To see (b), note that in this case $P_{\lambda,\mu}$ is a frame for $\cH$ $\Leftrightarrow$ $(g_\lambda,\Lambda)$ and $(g_\mu,\Lambda^{\circ})$ form a frame $\Leftrightarrow$ $(g_\lambda,\Lambda)$ and $(g_\mu,\Lambda^{\circ})$ form a Riesz basis $\Leftrightarrow$ $\{P_{\lambda,\mu} = g_\lambda^{\ast}\otimes  g_\mu$, $(\lambda,\mu) \in\Lambda\times\Lambda^{\circ}\}$ is a Riesz basis for $\cH$. Furthermore, by the classical Balian-Low theorem, if a Gabor system is an  $\mathbf{L}^2$-frame and at the same time a Riesz sequence (hence an  $\mathbf{L}^2$-Riesz basis), then the generating window $g$ cannot be in $S_0$\footnote{More precisely, $g$ cannot even be in the space of continuous functions in the Wiener space $W(\mathbb{R})$, see~\cite[Theorem~8.4.1]{gr01}}.
\hfill\foorp.

We may next ask, when the projection operators form a Riesz sequence, 
if  the reconstruction windows are TF-shifted
versions of a single window $h$ on the
adjoint lattice of $\Lambda =  b_0 \mathbb{Z}\times \nu_0 \mathbb{Z}$. In fact, in this case, 
the matrix $\Gamma$ turns out to enjoy quite a simple form.
To fix some notation, let 
\begin{eqnarray*}
\cA_{mn}(b,\nu) =& \sum_{k,\ell}e^{2i\pi m[\nu-\ell/\nu_0]}\ 
\overline{\cV_g h}(b-k/\nu_0,\nu-\ell/b_0)\\
\qquad\qquad\times & \ 
{\cV_g h}(b-(k-m)/\nu_0,\nu-(\ell-n)/b_0), 
\end{eqnarray*}
and introduce the right twisted convolution operator
$$
K_\cA^\natural(b,\nu): \cM(b,\nu)\to \cM(b,\nu)\natural \cA(b,\nu).
$$
\begin{corollary}
Let $g,h \in S_0$ as well as $b_0,\nu_0$ be given. Furthermore, let
$h^{(j)} = \pi(\frac{m}{\nu_0}, \frac{n}{b_0})h$. Then the
 variational equations read 
\begin{equation}
\label{eq:DefA}
\cM(b,\nu) \natural \cA(b,\nu) = \cB(b,\nu)\ . 
\end{equation} 
Hence, if  for all $b,\nu\in\RR^2$, the discrete right twisted
convolution operator $K_\cA^\natural$ is invertible, then the family $P_{\lambda,\mu} = g_\lambda^{\ast}\otimes  g_\mu$, $(\lambda,\mu) \in\Lambda\times\Lambda^{\circ}$ forms a Riesz sequence and the best MGM
approximation of an Hilbert-Schmidt operator with spreading function
$\eta$ is given by the family of transfer functions
$$
\cM_{mn}(b,\nu) = \left[(K_A^\natural(b,\nu))\inv \cB(b,\nu)\right]_{mn}\ ,
$$
where $\cB$ is given in~(\ref{fo:tototo}).
\end{corollary}
\underline{Proof:} As in the proof of Corollary~\ref{Co:ProdFram2}, we may derive the given form of $\Gamma$ by direct calculation, achieving the final form by noting that for $\Lambda= b_0\mathbb{Z}\times\nu_0\mathbb{Z}$, the adjoint lattice is given by $\Lambda^{\circ}=\frac{1}{\nu_0}\mathbb{Z}\times\frac{1}{b_0}\mathbb{Z}$.\hfill\foorp.\\

We close this section with some results of numerical experiments testing the approximation quality of MGMs for slowly time-varying systems.  \eqref{fo:glou}.

\begin{example}
We study the approximation of a (slowly) time-varying operator.
The operator has been generated by perturbing a time-invariant operator.
The spreading function is  shown in the upper display of Figure~\ref{FI2}.
The signal length is $32$, time- and frequency-parameters are $b_0 = 4$
and $\nu_0 = 4$, such that the redundancy of the  Gabor frame used in
the MGM approximation is $2$. The approximation is then realized in
several steps for two different schemes. Scheme~1 adds   three  synthesis window corresponding to a frequency-shift by $4$, a time-shift by  $4$ and a time-frequency-shift by $(4,4)$. The first step~$1$ calculates the
regular Gabor multiplier approximation. Step~$2$ adds one (only
frequency-shift) and so on. The rank of the resulting operator families is $64$, 
 $128$ for both step~2 and step~3 (adding either time- or frequency-shift) and  $256$ (time-shifted, frequency-shifted and time-frequency-shifted window added). The resulting approximation-errors are given by the solid line in the lower display of Figure~\ref{FI2}.\\
 Scheme~2 considers synthesis windows shifted in time and frequency on the sub-lattice generated by $a = 8, b = 8$, the resulting families having rank
 $64$,         $256$ and         $576$. Here, we only plot the results for the case corresponding to three and eight additional synthesis windows, respectively. The results are given by the dotted line.\\
  For comparison, an approximation with a regular Gabor multiplier with
redundancy $8$, i.e. an approximation family of rank $256$, has been
performed. The approximation error for this situation is the diamond in the middle of the display.
\begin{figure}
\centerline{\includegraphics[scale = .8]{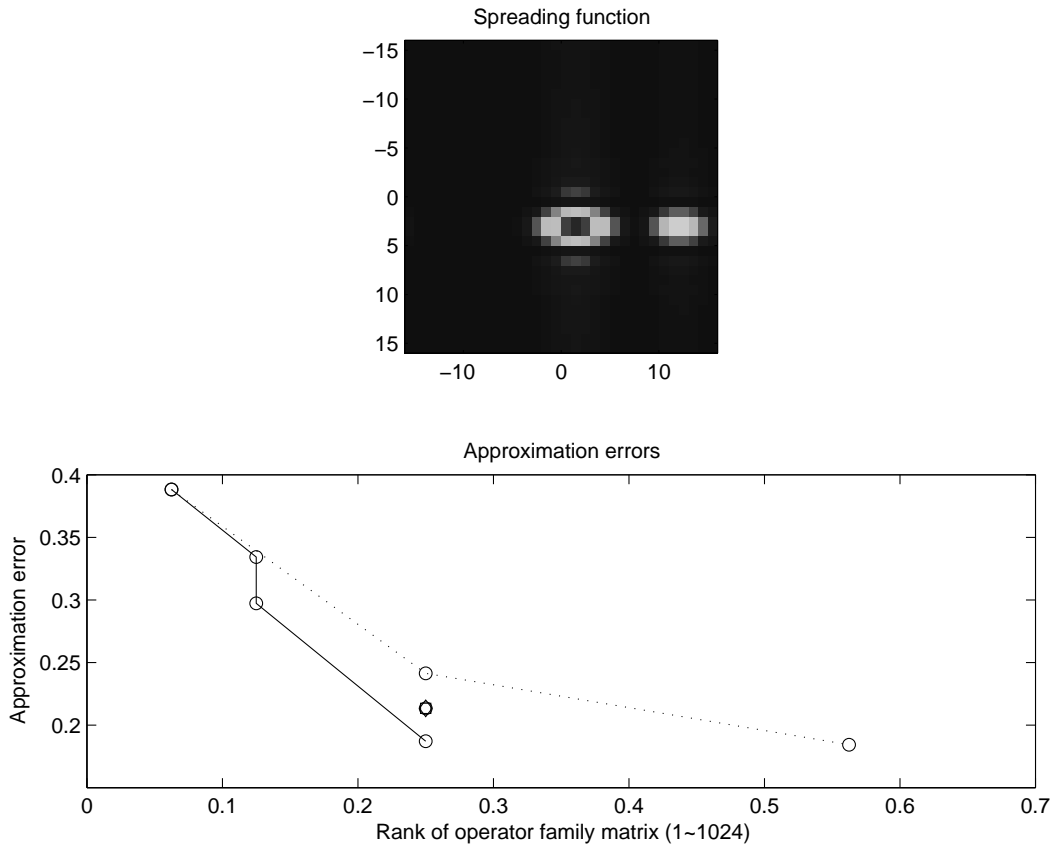}}
\caption{{\it Approximation by Multiple Gabor Multipliers }}
\label{FI2}
\end{figure}
 It is easy to see that,
depending on the behavior of the spreading function,  different schemes
perform advantageously for a certain redundancy. Note that for scheme~1, the best MGM with the same rank as the regular Gabor
multiplier performs better than the latter.  In the case of the
present operator, scheme~2 performs the "wrong" time-frequency
shifts on the synthesis windows in order to capture important
characteristics of the operator. However, in a different
setting, this scheme might be favorable (e.g. if an echo with a
longer delay is present).
\end{example}
 The example shows, that the choice of an appropriate sampling
scheme for the synthesis windows is extremely important in order to achieve a good and
efficient approximation by MGM. An optimal sampling scheme
depends on the analysis window's STFT, the lattice used in the
analysis and on the behavior of the operator's spreading function, which reflects the amount of delay and Doppler-shift created by the operator.
 Additionally, structural properties of the family of projections operators used in the approximation, based on the results in this section, have to be exploited to achieve numerical efficiency. An algorithm for optimization of these parameters
is currently under development.\footnote{P.~Balazs, M.~D\"orfler,
F.~Jaillet and B.~Torr\'esani. An optimized sampling scheme for
generalized Gabor multipliers. {\em {I}n preparation}, 2008.}

\subsection{Varying the synthesis window: TST spreading functions}
\label{Se:TSTSF}

We next turn to the special case of  separable  functions
$\bm_j (\lambda ) = m_1(\lambda)m_2(j)$ for the mask in the
definition \eqref{fo:multimulti} of MGMs.
In this case the resulting operator is of the form
\[
\MM f  = \sum_{\lambda} m_1(\lambda )\rho (\lambda )
(\sum_j m_2 (j) ( g^{\ast}\otimes  h^j))(f)
=\sum_{\lambda} m_1(\lambda )\rho (\lambda )\PP_m f\ ,
\]
where  $\PP_m f = \sum_j m_2 (j) \langle f, g\rangle  h^j$ and 
$\rho (\lambda )$  denotes a tensor product of
time-frequency shifts:
\[\rho (\lambda ) H := \pi(\lambda )H\pi^{\ast}(\lambda ).\]
Hence, the spreading function of $\MM$ is given by 
\[
\eta_{\MM} = \cM\cdot \eta_{\PP_m}, 
\]
where $\cM$ is the discrete  symplectic FT of $m_1$.
If the reconstruction windows are given by
$h^j = \pi (\mu^j )h$, $\mu_j = (b_j,\nu_j )$,  this becomes
\[
\eta_{\MM} = \cM\cdot\sum_j m_2(j) \cV_g h(\lambda -\mu_j )
e^{-2\pi i (\nu -\nu_j )b_j}\ .
\]
Motivated by this result, we introduce the following definition.
\begin{Def}[TST spreading functions]
Let $\phi$ be a given function from the function spaces $(S_0(\RR^2),
\mathbf{L}^2(\RR^2),S_0'(\RR^2))$
and let $b_1,\nu_1$ denote positive numbers.
Let $\mathbf{ \alpha}$ be in $\ell^1(\mathbb{Z}^2)$. A spreading
function $\eta=\eta_H$ of $H\in (\cBB,\cHH, \cBB')$,
that may be written as
\begin{equation}
\label{fo:TST}
\eta(b,\nu) = \sum_{k,\ell} \alpha_{k\ell} \phi(b-kb_1,\nu-\ell\nu_1)
e^{-2i\pi(\nu-\ell\nu_1)kb_1}\ 
\end{equation}
will be called {\em Twisted Spline Type} function (TST for
short).\end{Def} 
\begin{remark}\rm
By $\alpha$ in $\ell^1$, the series defining $\eta$ is
absolutely convergent in  $(S_0,\mathbf{L}^2,S_0' )$.
For $\ell^2$-sequences $\alpha$, we obtain an
$\mathbf{L}^2$-function $\eta$ for $\phi\in\mathbf{L}^2$.
\end{remark} 

TST functions are nothing but spline type functions (following the
terminology introduced in~\cite{Feichtinger02wavelet}), in which usual
(Euclidean) translations are replaced with the natural
(i.e. $\HH$-covariant) translations on the phase space $\PP$.
In fact, by writing
$\alpha (b,\nu ) = \sum_{k,\ell} \alpha_{k\ell} \delta(b-kb_1,\nu-\ell\nu_1)$,
the TST spreading function  may be written as a twisted convolution:
$\eta (b,\nu ) = \alpha\natural  \phi$. This
leads to the following property of operators associated with TST
spreading functions. 

\begin{lemma}
An operator $H\in (\cBB,\cHH, \cBB')$ possesses 
a TST spreading function
$\eta\in (S_0, \mathbf{L}^2,S_0')$
as in~(\ref{fo:TST}) if and only if it is of the form
\begin{equation}
H_{\eta} =  \sum_{k,\ell} \alpha_{k\ell} \pi(kb_1,\ell\nu_1) H_\phi\ ,
\end{equation}
where $H_\phi $ is the linear operator 
with spreading function $\phi$. 
\end{lemma}

\underline{Proof:} 
The proof consists of a straight-forward computation which may be
spared by noting that we have, by \eqref{fo:opTwist}:
\[
H_{\eta} =H_{\alpha\natural  \phi}= H_{\alpha}\cdot H_{\phi}
= \sum_{k,\ell} \alpha_{k\ell}  \pi(kb_1,\ell\nu_1) H_{\phi}\ .
\]
\hfill\foorp


As before, we are particularly interested in the situation of
the synthesis windows being given by time-frequency shifted versions
of a single window: $h^j = \pi(\mu_j ) h$.
In a next step we note, that under the condition
$\pi(\lambda ) \pi(\mu_j ) = \pi(\mu_j )\pi(\lambda )$,
i.e., $\mu_j\in\Lambda^{\circ}$, the MGM with separable multiplier
results in a TST spreading function with a Gabor multiplier as
basic operator $H_\phi$.
\begin{lemma}
\label{Le:MGM_TST}
Assume that a MGM $\mathbb{M}$ with multiplier
$m_j( \lambda ) = m_1 (\lambda ) m_2(j)$ is given. If the
synthesis windows $h^j$ are given by $h^j = \pi (\mu _j ) h$,
with $\mu_j\in\Lambda^{\circ}$, then 
\[
\mathbb{M} = \sum_j m_2 (j) \pi (\mu_j ) \mathbb{M}^G_{m_1;g,h},
\]
i.e., here, the operator $H_{\phi}$ is given by a regular Gabor
multiplier with mask $m_1$ and synthesis window $h$.
\end{lemma}
\begin{remark}\rm
Comparing the expression in the previous lemma to the expression
$\mathbb{M} = \sum_{\lambda} m_1(\lambda )\rho (\lambda )\PP_m $ for
the same operator, we note that in this situation, the operator may
either be interpreted as a (weighted) sum of Gabor multipliers or
as a Gabor multiplier with a generalized projection operator
$\mathbb{P}_{\bm}$ in the   synthesis process. In this situation,
we may ask, whether the family of generalized projection operators,
$\{\rho (\lambda )\mathbb{P}_{\bm}\}_{\lambda\in\Lambda}$ form a
frame or Riesz basis for their linear span. In fact, if $\bm$ is
in $\ell^1$ and $g,h\in S_0$, this question is easily answered by
generalizing the result proved in~\cite[Theorem~3.2]{bepf06}.
Here, $\{\rho (\lambda )\mathbb{P}_{\bm}\}_{\lambda\in\Lambda}$ is
either a Riesz basis or not a frame for its closed linear span.
Furthermore, there exists $r>0$ such that
$\{\rho (\alpha\lambda )\mathbb{P}_{\bm}\}_{\lambda\in\Lambda}$ is
a Riesz basis for its closed linear span whenever $\alpha>r$.
\end{remark}

In generalizing the result of Lemma~\ref{Le:MGM_TST}, it  is a natural next step to assume that the basic function  $\phi$ entering
in the composition of $\eta$ is the spreading function of a Gabor
multiplier (at least in an approximate sense). According to
the discussion of Section~\ref{se:mult}, this essentially means
that $\phi$ is sufficiently well concentrated in the
time-frequency domain. \\
(In the sequel we will write $\pi_{mn}$ for $\pi (mb_0,n\nu_0)$
whenever the applicable lattice constants are sufficiently clear
from the context.)\\
Hence, we assume that a Gabor multiplier $H_\phi$, as defined
in~\eqref{fo:GMdef} is given. We may formally compute
\begin{eqnarray}\label{eq:TSTGabmul}
Hf &=& \sum_{k,\ell} \alpha_{k\ell} \pi_{k\ell}\sum_{m,n} \bm(m,n)
\cV_gf(mb_0,n\nu_0) \pi (mb_0,n\nu_0) h\\\nonumber
&=& \sum_{m,n} \bm(m,n)\cV_gf(mb_0,n\nu_0)
\sum_{k,\ell} \alpha_{k\ell}\pi(kb_1,\ell\nu_1)\pi (mb_0,n\nu_0) h
\end{eqnarray}
Based on this expression, one may pursue two different choices
of the sampling-points $(kb_1, \ell\nu_1)$. First, in extension
of the result given in Lemma~\ref{Le:MGM_TST}, we assume that
the sampling points are associated to the adjoint lattice
$\Lambda^{\circ} = \frac{1}{\nu_0}\mathbb{Z}\times  \frac{1}{b_0}\mathbb{Z}$
of $\Lambda = b_0\mathbb{Z}\times \nu_0\mathbb{Z}$.
The second choice of sampling points on the original  lattice leads
to a construction as introduced in~\cite{Dorfler07spreading} as
Gabor twisters and will not be further discussed in the present contribution.\\

The following theorem extends the result given in
Lemma~\ref{Le:MGM_TST} to the case in which the sampling
points in the TST expansion are chosen from a lattice
containing the adjoint lattice. It turns out that the TST
spreading function then leads to a representation as a
sum of Gabor multipliers.

\begin{theorem}\label{Th:TSTfacts}
Let $b_0,\nu_0\in\RR^+$ generate the time-frequency lattice $\Lambda$, and
let $\Lambda^\circ$ denote the adjoint lattice.
Let $g,h\in S_0(\RR)$ denote respectively Gabor
analysis and synthesis windows,
such that the $\cU$ condition~(\ref{fo:Ucond}) is fulfilled.
Let $H$ denote the operator in $(\cBB, \cHH, \cBB')$ defined
by the twisted spline type spreading function $\eta$ as
in~(\ref{fo:TST}), with $b_1,\nu_1\in\RR^+$.
\begin{enumerate}
\item
Assume that $b_1$ and $\nu_1$ are multiple of the dual lattice constants.
Then $H$ is a Gabor multiplier, with analysis window $g$,
synthesis window
\begin{equation}
\gamma =\sum_{k,\ell} \alpha_{k\ell}\pi(kb_1,\ell\nu_1) h\ ,
\end{equation}
and transfer function
\begin{equation}
\label{fo:TST.transfer.function}
\bm(m,n) = b_0\nu_0\int_\square \cM(b,\nu) e^{-2i\pi(n\nu_0b - mb_0\nu)}\, db
d\nu\ ,
\end{equation}
with $\square$ the fundamental domain of the adjoint
lattice $\Lambda^\circ$, and
\begin{equation}
\label{fo:TST.transfer.function2}
\cM(b,\nu) = \frac{\sum_{k,\ell=-\infty}^\infty
\overline{\cV_gh}\left(b+k/\nu_0,\nu + \ell/b_0\right)
\phi\left(b+k/\nu_0,\nu + \ell/b_0\right)}
{\sum_{k,\ell=-\infty}^\infty
\left|\cV_gh\left(b+k/\nu_0,\nu + \ell/b_0\right)\right|^2}
\end{equation}
\item
Assume that the lattice generated by $b_1$ and $\nu_1$ contains the
adjoint lattice:
\begin{equation}
\label{fo:super.dual.lattice}
b_1 = \frac1{p\nu_0}\ ,\qquad \nu_1 = \frac1{qb_0}\ .
\end{equation}
Then $H$ may be written as a finite sum of Gabor multipliers
\begin{equation}
\label{fo:super.lattice.approx}
Hf = \sum_{i=1}^p \sum_{j=1}^q \bigg(
\sum_{m\equiv i\,[{ \rm mod}\,p]}\ \sum_{n\equiv j\,[{ \rm mod}\,q]}
\bm(m,n)\cV_gf(mb_0,n\nu_0) \pi_{mn}\bigg) \gamma_{ij}\ ,
\end{equation}
with at most $p\cdot q$ different synthesis windows $\gamma_{ij}$
and the transfer function given in~(\ref{fo:TST.transfer.function})
and~(\ref{fo:TST.transfer.function2}).
\end{enumerate}
\end{theorem}

\underline{Proof:} \\
Let us formally compute
\begin{eqnarray*}
Hf &=& \sum_{m,n} \bm(m,n)\cV_gf(mb_0,n\nu_0)
\sum_{k,\ell} \alpha_{k\ell}\pi(kb_1,\ell\nu_1)\pi_{mn} h\\
&=& \sum_{m,n} \bm(m,n)\cV_gf(mb_0,n\nu_0) \pi_{mn} \gamma_{mn}\ ,
\end{eqnarray*}
where
\begin{equation}
\nonumber
\gamma_{mn} = \sum_{k,\ell} \alpha_{k\ell}
e^{2i\pi [knb_0\nu_1 - \ell m \nu_0 b_1]} \pi(kb_1,\ell\nu_1) h
\end{equation}
Now observe that if $(b_1,\nu_1)\in\Lambda^\circ$,
one obviously has
$$
\gamma_{mn} = \sum_{k,\ell} \alpha_{k\ell}\pi(kb_1,\ell\nu_1) h =
\gamma_{00}\ , \mbox{ for } (m,n)\in\ZZ^2\ ,
$$
i.e. the above expression for $Hf$ involves a single synthesis window
$\gamma=\gamma_{00}$. Therefore, in this case, $H$ takes the form of a
standard Gabor multiplier, with fixed time-frequency transfer
function, and a synthesis window prespribed by the coefficients in the
TST expansion. This proves the first part of the theorem.

Let us now assume that the TST expansion of the spreading function is
finer than the one prescribed by the lattice $\Lambda^\circ$, but
nevertheless the lattice $\Lambda_1 = \ZZ b_1\times\ZZ\nu_1$ contains
$\Lambda^\circ$. In other words, there exist positive integers $p,q$
such that \eqref{fo:super.dual.lattice} holds.\\
We then have
\begin{equation}
\gamma_{mn} = \sum_{k,\ell} \alpha_{k\ell}
e^{2i\pi [\frac{knp-lmq}{pq} ]} \pi(kb_1,\ell\nu_1) h
\end{equation}
and it is readily seen that there are at most $pq$
different synthesis windows $\gamma_{ij}$,
\begin{equation}
\label{fo:super.synth.win}
  \gamma_{ij} = \gamma_{m\,[{\rm mod}\,p],n\,[{\rm mod}\,q]}\ , 
i = 1,\ldots , p; \ j = 1,\ldots , q .
\end{equation}
The operator $H$ may hence be written as a sum of Gabor multipliers, 
with one prescribed time-frequency transfer function, which is
sub-sampled on several sub-lattices of the lattice $\Lambda$:
\[
\Lambda_{ij} = (p b_0\cdot\mathbb{Z}+i\cdot b_0)\times 
(q \nu_0\cdot\mathbb{Z}+j\cdot \nu_0),\, 
i = 0,\ldots, p-1; j = 0, \ldots , q-1, \]
and a single synthesis window per sub-lattice as given
in~\eqref{fo:super.synth.win}.
The resulting expression for $H$ is hence as given
in~\eqref{fo:super.lattice.approx}. \\

The expression for the transfer function is derived in analogy to the
case discussed in Section~\ref{se:mult}.
\hfill\foorp.

\begin{remark}\rm
Let us observe that in this approximation, the time-frequency transfer
function $\bm$ is completely characterized by the function $\phi$ used
in the TST expansion. The choice of $\phi$ therefore imposes a fixed
mask for the multipliers that come into play in
equation~(\ref{fo:super.lattice.approx}).

\end{remark}

\begin{example}
We first assume, that for a given primal lattice $\Lambda =
b_0\mathbb{Z}\times \nu_0\mathbb{Z}$, the representation of
a spreading function $\eta$ is given 
by $5$ building blocks:
\[\eta (b,\nu ) =\sum_{k = -1}^1 \alpha_{k0} \phi (b
-\frac{k}{\nu_0},\nu )+\sum_{\ell = -1}^1 \alpha_{0\ell}    \phi (b ,\nu
-\frac{\ell}{b_0}).\] 
In this case, we obtain a single Gabor multiplier with synthesis window 
\[
\gamma_{00} = \sum_{k = -1}^1 \alpha_{k0} \pi (\frac{k}{\nu_0},0)h
+\sum_{\ell = -1}^1 \alpha_{0\ell} \pi (0, \frac{\ell}{b_0})h\ .
\]
{
If we add the windows $ \phi (b \pm\frac{1}{2\nu_0},\nu
\pm\frac{1}{2b_0})$ to the representation of $\eta$, we are now
dealing with the finer lattice $\Lambda =
\frac{1}{2\nu_0}\mathbb{Z}\times\frac{1}{2b_0} \mathbb{Z}$ and we
obtain the sum of $4$ Gabor multipliers with the following synthesis windows: 
\begin{eqnarray*}
\gamma_{00} &=& \sum_{k = -1}^1 \alpha_{k0}
\pi\! \left(\frac{k}{2\nu_0},0\right)h
+\sum_{\ell = -1}^1 \alpha_{0\ell}
\pi\! \left(0, \frac{\ell}{2b_0}\right)h, \\
\gamma_{01} &= &\sum_{k = -1}^1 \alpha_{k0}e^{\pi ik}
\pi\!\left(\frac{k}{2\nu_0},0\right)h +
\sum_{\ell = -1}^1 \alpha_{0\ell}
\pi\! \left(0,\frac{\ell}{2b_0}\right)h, \\
\gamma_{10} &=& \sum_{k = -1}^1 \alpha_{k0}
\pi\!\left(\frac{k}{2\nu_0},0\right)h
+\sum_{\ell = -1}^1e^{\pi i\ell}\alpha_{0\ell}
\pi\! \left(0, \frac{\ell}{2b_0}\right)h\ ,\\
\gamma_{11} &=& \sum_{k = -1}^1 \alpha_{k0}e^{\pi ik}
\pi\!\left(\frac{k}{2\nu_0},0\right)h +
\sum_{\ell = -1}^1 \alpha_{0\ell}e^{-\pi i\ell}
\pi\! \left(0,\frac{\ell}{2b_0}\right)h\ ,
\end{eqnarray*}
and corresponding lattices: 
$\Lambda_{00} = 2\mathbb{Z}b_0\times 2\mathbb{Z}\nu_0$, 
$\Lambda_{01} = 2\mathbb{Z}b_0\times (2\mathbb{Z}+1)\nu_0$,
$\Lambda_{10} = (2\mathbb{Z} +1)b_0\times 2\mathbb{Z}\nu_0$, and
$\Lambda_{11} = (2\mathbb{Z}+1)b_0\times (2\mathbb{Z}+1)\nu_0$.
}\end{example}

\medskip

It is important to note, that in both cases described in
Theorem~\ref{Th:TSTfacts} as well as the above example,  the transfer
function $\bm$ can be calculated as the best approximation by  a
regular Gabor multiplier  - a procedure which may be efficiently realized using \eqref{Eq:MBestGM}. Fast algorithms for this 
exist in the literature, see~\cite{fehakr04}, however, the method derived in Section~\ref{Se:tfmult} appears to be faster. 
\section{Conclusions and Perspectives} 
Starting from an operator representation in the continuous time-frequency domain via a twisted convolution,  we have introduced generalizations of conventional time-frequency multipliers in order to overcome the restrictions of this model in the approximation of general operators. 
The model of multiple Gabor multipliers in principle allows the representation of any given linear operator. However, in oder to achieve computational efficiency as well as insight in the operator's characteristics, the parameters used in the model must be carefully chosen. An algorithm choosing the optimal sampling points for the family of synthesis windows, based on the spreading function, is the topic of ongoing research. On the other hand, the model of twisted spline type functions allows the approximation of a given spreading function and results in an adapted window or family of windows. By refining the sampling lattice in the TST approximation, a rather wide class  of operators should be well-represented. The practicality of this approach has to be shown in the context of operators of practical relevance. All the results given in this work will also be applied  in the context of estimation rather than approximation.\\
As a further step of generalization, frame types other than Gabor frames may be considered. Surprisingly little is known about wavelet frame multipliers, hence it will be interesting to  generalize the achieved results  to the affine group. 
\bibliographystyle{abbrv}

\end{document}